\documentclass{amsart}
\usepackage{amssymb,amscd}

\newcounter{intro}

\newtheorem{theo}[intro]{Theorem}

\newtheorem{thm}{Theorem}[section]
\newtheorem{lem}[thm]{Lemma}
\newtheorem{prop}[thm]{Proposition}
\newtheorem{cor}[thm]{Corollary}

\theoremstyle{remark}
\newtheorem{rem}[thm]{Remark}
\newtheorem*{merci}{Acknowledgements}

\numberwithin{equation}{section}

\newcommand{\cref}[1]{Corollary~\ref{#1}}

\newcommand{\lref}[1]{Lemma~\ref{#1}}
\newcommand{\pref}[1]{Proposition~\ref{#1}}

\newcommand{\tref}[1]{Theorem~\ref{#1}}
\newcommand{\trefs}[1]{Theorems~\ref{#1}}
\newcommand{\sref}[1]{Section~\ref{#1}}

\newcommand{\Hyp}{\mathbb{H}}
\newcommand{\R}{\mathbb{R}}
\newcommand{\C}{\mathbb{C}}
\renewcommand{\H}{\mathbb{H}}
\newcommand{\K}{\mathbb{K}}
\renewcommand{\O}{\mathbb{O}}
\newcommand{\N}{\mathbb{N}}
\newcommand{\Z}{\mathbb{Z}}
\newcommand{\bS}{\mathbb{S}}

\newcommand{\boch}{\mathcal{B}}
\newcommand{\cour}{\mathcal{R}}
\newcommand{\cH}{\mathcal{H}}

\newcommand{\cN}{\mathcal{N}}

\newcommand{\Ag}{\mathfrak{a}}
\newcommand{\Gg}{\mathfrak{g}}
\newcommand{\Mg}{\mathfrak{m}}
\newcommand{\Pg}{\mathfrak{p}}
\newcommand{\Lg}{\mathfrak{l}}
\newcommand{\Ng}{\mathfrak{n}}

\newcommand{\Kg}{\mathfrak{k}}

\let\al=\alpha

\let\de=\delta
\let\De=\Delta
\let\eps=\varepsilon
\let\ga=\gamma
\let\Ga=\Gamma
\let\la=\lambda
\let\Om=\Omega
\let\om=\omega

\let\si=\sigma
\let\Si=\Sigma

\let\priv=\smallsetminus
\let\wh=\widehat

\DeclareMathOperator{\End}{End}
\DeclareMathOperator{\Hom}{Hom}
\DeclareMathOperator{\Id}{Id}
\DeclareMathOperator{\im}{Im}
\DeclareMathOperator{\Ind}{Ind}
\DeclareMathOperator{\rank}{rank}
\DeclareMathOperator{\re}{Re}
\DeclareMathOperator{\spec}{spec}
\DeclareMathOperator{\supp}{supp}

\DeclareMathOperator{\vol}{vol}

\def\cro#1#2{\mathrel{\langle {#1},{#2}\rangle}}
\def\norm#1{\Vert {#1}\Vert}

\begin{document}

\title[On the differential form spectrum of hyperbolic manifolds]
{On the differential form spectrum\\ of hyperbolic manifolds}

\author{Gilles Carron}
\address{Laboratoire de Math\'ematiques Jean Leray (UMR 6629), Universit\'e de Nantes, 
2, rue de la Houssini\`ere, B.P.~92208, 44322 Nantes Cedex~3, France}
\email{Gilles.Carron@math.univ-nantes.fr}
\author{Emmanuel Pedon}
\address{Laboratoire de Math\'ematiques (UMR 6056), Universit\'e de Reims, Moulin de la
Housse, B.P.~1039, 51687 Reims Cedex~2, France}
\email{emmanuel.pedon@univ-reims.fr}

\subjclass[2000]{Primary 53C35, 58J50; Secondary 22E40, 34L15, 57T15}
\keywords{Hyperbolic spaces, hyperbolic manifolds, locally symmetric spaces, cohomology,
Hodge-de~Rham Laplacian, spectral theory}

\date{September 15, 2004}

\begin{abstract}
We give a lower bound for the bottom of the $L^2$ differential form spectrum on hyperbolic
manifolds, generalizing thus a well-known result due to Sullivan and Corlette in the function
case. Our method is based on the study of the resolvent associated with the Hodge-de~Rham
Laplacian and leads to applications for the (co)homology and topology of certain
classes of hyperbolic manifolds. 
\end{abstract}

\maketitle

\section{Introduction}

Let $G/K$ be a Riemannian symmetric space of noncompact type, and let $\Ga$ be a discrete,
torsion-free subgroup of $G$. Thus $\Ga\backslash G/K$ is a locally Riemannian symmetric space
with nonpositive sectional curvature. Most of this article concerns the rank one case, i.e. when
$G/K$ is one of the hyperbolic spaces $\Hyp^n_\R$, $\Hyp^n_\C$, $\Hyp^n_\H$ or $\Hyp^2_\O$.
In that situation, the quotients $\Ga\backslash\Hyp^n_\K$ are usually called hyperbolic
manifolds, and we normalize the Riemannian metric so that the corresponding pinched sectional
curvature lies inside the interval $[-4,-1]$.

We denote by $2\rho$ the exponential rate of the volume growth in
$\Hyp^n_\K$:
$$2\rho=\lim_{R\to +\infty} \frac{\log \vol B(x,R)}{R},$$
and let $\de(\Ga)$ be the critical exponent of the Poincar\'e series associated with $\Ga$, i.e.
$$\delta(\Gamma)=\inf\{s\in\R \textrm{\ such\ that\ } \sum_{\gamma\in\Gamma} e^{-s
d(x,\gamma y)} <+\infty\},$$
where $(x,y)$ is any pair of points in $\Hyp^n_\K$ and $d(x,\ga y)$ is the geodesic distance
from $x$ to $\ga y$. It is well-known that $0\le\de(\Ga)\le 2\rho$.

For any (locally) symmetric space $X$ considered above, let $\lambda_0^p(X)$ be the bottom
of the $L^2$ spectrum of the Hodge-de~Rham Laplacian $\De_p$ acting on compactly supported
smooth differential $p$-forms of $X$. In other words,
$$\lambda_0^p(X)=\inf_{u\in C^\infty_0(\wedge^p T^*X)}
\frac{(\De_p u|u)_{L^2}} {\|u\|_{L^2}^2}.$$

Let us recall the following beautiful result, due to D.~Sullivan (\cite{sullivanl},
Theorem~2.17) in the real case and to K.~Corlette (\cite{corlette},
Theorem~4.2) in the remaining cases (see also \cite{Elstrodt}, \cite{patterson} and
\cite{CdV} in the case of $\Hyp^2_\R$):

\begin{theo}\label{SC}
\begin{enumerate}
\item If $\delta(\Gamma)\le \rho$, then $\lambda_0^0(\Gamma \backslash\Hyp_\K^n)=\rho^2$.
\item If $\delta(\Gamma)\ge \rho$, then 
$\lambda_0^0(\Gamma \backslash\Hyp_\K^n)=\delta(\Gamma)\left(2\rho-\delta(\Gamma)\right)$.
\end{enumerate}
\end{theo}

The main goal of our paper is to extend this result to the case of differential forms, although
we are aware that getting such a simple statement is hopeless. 
For instance, when $\Ga$ is cocompact 
the zero eigenspace of the Hodge-de~Rham Laplacian $\Delta_p$ acting on
$\Gamma \backslash\Hyp_\K^n$
is isomorphic to the $p$-th cohomology group of $\Gamma \backslash\Hyp_\K^n$, and contains
therefore some information on the topology of this manifold. Thus one does not
expect to compute the bottom of the spectrum of $\Delta_p$ only in terms of the critical exponent,
since we always have $\de(\Ga)=2\rho$ in the cocompact case.

Nevertheless, we are able to give lower bounds for $\lambda_0^p(\Gamma
\backslash\Hyp_\K^n)$. In order to state our first
result, we set 
$d=\dim_\R(\K)$ and denote by $\al_p$ the bottom of the continuous $L^2$ spectrum of $\De_p$ on the
hyperbolic space $\Hyp^n_\K$.

\begin{theo} \label{princ}
\begin{enumerate}
\item Assume that $p\not=\frac{dn}{2}$.
\begin{enumerate}
\item If $\delta(\Gamma)\le \rho$, then
$\lambda_0^p(\Gamma \backslash\Hyp_\K^n)\ge \alpha_p$.
\item If $\rho\le\delta(\Gamma)\le \rho+\sqrt{\al_p}$, then
$\lambda_0^p(\Gamma \backslash\Hyp_\K^n)\ge
\alpha_p-\left(\delta(\Gamma)-\rho\right)^2$.
\end{enumerate}
\item Assume that $p=\frac{dn}{2}$.
\begin{enumerate}
\item If $\delta(\Gamma)\le \rho$, then either $\lambda_0^p(\Gamma \backslash\Hyp_\K^n)=0$ or
$\lambda_0^p(\Gamma \backslash\Hyp_\K^n)\ge \alpha_p$.
\item If $\rho\le\delta(\Gamma)\le \rho+\sqrt{\al_p}$, then 
either $\lambda_0^p(\Gamma \backslash\Hyp_\K^n)=0$ or
$\lambda_0^p(\Gamma \backslash\Hyp_\K^n)\ge
\alpha_p-\left(\delta(\Gamma)-\rho\right)^2$.
\end{enumerate}
Moreover, if $\delta(\Gamma)<\rho+\sqrt{\al_p}$ the possible
eigenvalue $0$ is discrete and spectrally isolated.
\end{enumerate}
\end{theo}

When $\delta(\Gamma)>\rho+\sqrt{\al_p}$, assertions (b) are still valid, but yield a triviality
since the spectrum must be non negative.

U.~Bunke and M.~Olbrich pointed out to us that, in the case of convex cocompact
subgroups $\Ga$, \tref{princ} could be obtained as a
consequence of Theorem~4.7 in \cite{BOJFA} or Theorem~1.8 in \cite{BO}.
However, besides it works in any case, our proof follows a completely different path, relying on
an estimate for the resolvent associated with $\De_p$ on $\Hyp^n_\K$. 
In particular, we are also
able to discuss the nature of the continuous spectrum of $\De_p$ on $\Ga\backslash\Hyp^n_\K$
when $\de(\Ga)<\rho$ (see \pref{natspec}).

Considering the following large class of examples, we see that our estimates in \tref{princ} are
sharp when $\de(\Ga)\le\rho$.

\begin{theo} \label{sharp}
If $\delta(\Gamma)\le \rho$ and if the injectivity radius of
$\Gamma\backslash\Hyp_\K^n$ is not bounded (for instance if the limit set $\Lambda(\Ga)$ of
$\Ga$ is not the whole sphere at infinity $\bS^{dn-1}$), then 
$$\spec \left(\Delta_p,\Gamma\backslash\Hyp_\K^n\right)
=\spec \left(\Delta_p,\Hyp_\K^n\right)
=\begin{cases}
[\alpha_p,+\infty)&\text{if }p\not=\frac{dn}{2},\\
\{0\}\cup[\alpha_p,+\infty)&\text{if }p=\frac{dn}{2}.
\end{cases}$$
\end{theo}

Since the exact value of $\al_p$ is known except in the case of $\Hyp^2_\O$ 
(see \tref{valalphap}), \tref{princ} provides an explicit
vanishing result for the space of $L^2$ harmonic forms, from which we shall obtain several
corollaries, most of them having a topological flavour. For instance, we give sufficient
conditions for a hyperbolic manifold to have only one end (actually we also deal with general 
locally symmetric spaces whose isometry group satisfies Kazhdan's property). Denote as usual by
$H_p(\Gamma\backslash\Hyp_\K^n,\Z)$ the $p$-th homology space of
$\Gamma\backslash\Hyp_\K^n$ with coefficients in $\Z$.

\begin{theo} \label{bout1}
Let $\Gamma$ be a discrete and torsion-free subgroup of the isometry group
of a quaternionic hyperbolic space $\Hyp^n_\H$ or of the octonionic hyperbolic plane $\Hyp^2_\O$.
If all unbounded connected components of the complement of any compact subset of 
$\Gamma\backslash\Hyp_\K^n$ have infinite volume,
then $\Gamma\backslash\Hyp_\K^n$ has only one end, and 
$$H_{dn-1}(\Gamma\backslash\Hyp_\K^n,\Z)=\{0\}.$$
\end{theo}

\begin{theo}\label{bout2}
Let $\Gamma$ be a discrete and torsion-free subgroup of $SU(n,1)$, with $n\ge 2$. Assume that the
limit set
$\Lambda(\Gamma)$ is not the whole sphere at infinity
$\bS^{2n-1}$, that $\delta(\Gamma)<2n$, and that the injectivity radius of
$\Gamma\backslash\Hyp_\C^n$ has a positive lower bound. Then
$\Gamma\backslash\Hyp_\C^n$ has only one end, and
$$H_{2n-1}(\Gamma\backslash\Hyp_\C^n,\Z)=\{0\}.$$
\end{theo}

The first of these two theorems extends a previous result of K.~Corlette (\cite{corlette},
Theorem~7.1) in the convex cocompact setting. The second enables us to complement a rigidity result
due to Y.~Shalom (\cite{shalom}, Theorem~1.6; see also \cite{BCGpre}):

\begin{theo} \label{shalfin}
Assume that $\Gamma=A*_C B$ is a cocompact subgroup of $SU(n,1)$ (with $n\ge 2$) which is a free
product of subgroups $A$ and $B$ over an amalgamated subgroup $C$. Then either
$2n-1\le\delta(C)<2n$ and $\Lambda(C)=\bS^{2n-1}$, or $\delta(C)=2n$.
\end{theo}

Our article is organized as follows. Section~2 contains most of the notation and background
material that will be used in this article, and especially a fairly detailed introduction to $L^2$
harmonic analysis on the differential form bundle over hyperbolic spaces, from the representation
theory viewpoint, since this approach is the touchstone of our work. We also briefly comment on
the generalization of
\tref{SC} to general nonpositively curved locally symmetric spaces and quotients of Damek-Ricci
spaces.

Section~3 is devoted to the analysis of the resolvent 
$$R_p(s)=(\De_p-\al_p+s^2)^{-1},\quad\textrm{for }\re s>0,$$
associated with the Hodge-de~Rham Laplacian
on hyperbolic spaces. More precisely, we obtain a meromorphic continuation on a suitable
ramified cover of $\C$, prove estimates at infinity, and discuss the possible location of the poles
on the imaginary axis $\re s=0$ of $\C$.

In Section~4 we prove the spectral results announced above (\trefs{princ}, \ref{sharp}), and apply
them to derive several vanishing results for the cohomology. We also verify that our results
on the bottom of the spectrum are strictly better than the ones given by the Bochner-Weitzenb\"ock
formula and the Kato inequality.

Lastly, Section~5 contains the proof of all results dealing with the number of ends and the
homology of locally symmetric spaces, in particular of \trefs{bout1}, \ref{bout2}, \ref{shalfin}. 

Numerous comments and references will be given throughout the text.

\begin{merci} We are particularly grateful to M.~Olbrich for his careful reading of this
article and for the numerous comments he made on it. We would like also to thank G.~Besson,
G.~Courtois and S.~Gallot for communicating to us quite soon a result of
\cite{BCGpre}, as well as J.-Ph.~Anker, U.~Bunke, P.-Y.~Gaillard, E.~Ghys, L.~Guillop\'e and
F.~Laudenbach for useful remarks and fruitful discussions.
\end{merci}

\section{Notations and background material}

In this section, we shall collect some notations, definitions and preliminary facts which will
be used throughout the article. Although some of our results concern
general locally symmetric spaces of noncompact type, our paper essentially deals with (quotients
of) hyperbolic spaces, and we prefer therefore to restrict the following comprehensive
presentation to that case. Most of unreferred material can be found for instance in the classical
books
\cite{Helgason} and \cite{Knapp86}.

\subsection{Hyperbolic spaces}
\label{notation}

For $n\geq 2$ and $\K=\R,\C,\H$ or for $n=2$ and $\K=\O$, let $\Hyp^n_\K$ be the Riemannian
hyperbolic space of dimension $n$ over
$\K$. Recall that $\Hyp^n_\K$ is realized as the noncompact symmetric space
of rank one
$G/K$, where $G$ is a
connected noncompact semisimple real Lie group with finite centre (namely, the identity
component of the group of isometries of $\Hyp^n_\K$) and
$K$ is a maximal compact subgroup of $G$ which consists of elements fixed by a Cartan involution
$\theta$. More precisely,
$$\begin{array}{cccccc} {\rm if} & \K=\R & {\rm\ then\ }& 
G=SO_e(n,1)&{\rm\ and\ } &
K=SO(n);\\ {\rm if} & \K=\C & {\rm\ then\ } & G=SU(n,1)& {\rm\ 
and\ } & K=S(U(n)\times
U(1));\\ {\rm if} & \K=\H & {\rm\ then\ } & G=Sp(n,1)& {\rm\ 
and\ } & K=Sp(n)\times
Sp(1);\\ {\rm if} & \K=\O & {\rm\ then\ } n=2 {\rm\ and\ }& 
G=F_{4(-20)}& {\rm\  and\ } &
K=Spin(9).\\
\end{array}$$ (Other pairs $(G,K)$ may be taken to give the same 
quotient $G/K$.)

Let us begin with some algebraic structure of the Lie groups 
involved.  Let $\Gg$ and $\Kg$ be
the Lie algebras of
$G$ and $K$, respectively, and write
\begin{equation}\label{Cartang}
\Gg=\Kg\oplus\Pg
\end{equation}  for the Cartan decomposition of $\Gg$ (i.e. the 
decomposition of $\Gg$ into
eigenspaces for the eigenvalues $+1$, $-1$, respectively, of the Cartan involution $\theta$).
Recall that the subspace
$\Pg$ is thus identified with the tangent space 
$T_{eK}\left(G/K\right)\simeq\R^{dn}$ of
$\Hyp^n_\K=G/K$ at the origin, where $d=\dim_\R(\K)$.

Let $\Ag$ be a maximal abelian subspace of $\Pg$ ($\Ag\simeq\R$ since 
${\rm rank}(G/K)=1$),
with corresponding analytic Lie subgroup $A=\exp(\Ag)$ of $G$. Let 
$R(\Gg,\Ag)$ be the
restricted root system of the pair $(\Gg,\Ag)$, with positive 
subsystem $R^+(\Gg,\Ag)$
corresponding to the positive Weyl chamber
$\Ag_+\simeq (0,+\infty)$ in $\Ag$. It is standard that there exists 
a linear functional
$\al\in\Ag^*$ such that
\begin{eqnarray} R(\Gg,\Ag)&=&
\begin{cases}\{\pm\al\}&\textrm{if }\K=\R,\\ \{\pm\al,\pm 
2\al\}&\textrm{if }\K=\C,\H,\O,\\
\end{cases}\\
\textrm{and }R^+(\Gg,\Ag)&=&
\begin{cases}\{\al\}&\textrm{if }\K=\R,\\ \{\al,2\al\}&\textrm{if 
}\K=\C,\H,\O.\\
\end{cases}
\end{eqnarray} As usual, we write
$\Ng$ for the direct sum of positive root subspaces, i.e.
\begin{equation}\label{defn}
\Ng=
\begin{cases}\Gg_\al&\textrm{if }\K=\R,\\ 
\Gg_\al\oplus\Gg_{2\al}&\textrm{if }\K=\C,\H,\O,\\
\end{cases}
\end{equation} so that $\Gg=\Kg\oplus\Ag\oplus\Ng$ is an Iwasawa 
decomposition for $\Gg$. We
let also $N=\exp(\Ng)$ and $\rho=\frac{1}{2}(m_\al
\al+m_{2\al}2\al)$, where $m_\al=\dim_\R\Gg_\al=d(n-1)>0$ and $m_{2\al}=\dim_\R
\Gg_{2\al}=d-1\ge 0$. In the sequel, we shall use systematically the 
identification
\begin{equation}\label{identaR}
\begin{split}
\Ag^*_\C&\simeq\C,\\
\la\al&\mapsto\la.
\end{split}
\end{equation} 
In particular, we shall view $\rho$ as a real number, namely
\begin{equation}\label{rho}
\rho=\frac{d(n-1)}{2}+d-1=
\begin{cases}
\frac{n-1}{2}&\textrm{if }\K=\R,\\ n&\textrm{if }\K=\C,\\ 
2n+1&\textrm{if }\K=\H,\\
11&\textrm{if }\K=\O\textrm{ and }n=2.\\
\end{cases}
\end{equation} 
This number has also a well-known geometrical interpretation: if $h$ denotes the exponential rate
of the volume growth in $\Hyp^n_\K$, i.e. if
\begin{equation*}
h=\lim_{R\to \infty} \frac{\log \vol B(x,R)}{R},
\end{equation*}
(this quantity does not depend on $x\in \Hyp^n_\K$) then $h=2\rho$.

Next, let $H_0\in\Ag_+$ be such that $\al(H_0)=1$.  We define a symmetric 
bilinear form on $\Gg$ by
\begin{equation}\label{bilin}
\cro{X}{Y}=\frac{1}{B(H_0,H_0)}B(X,Y)=\frac{1}{2(m_\al+4m_{2\al})}B(X, 
Y),
\end{equation} where $B$ is the Killing form on $\Gg$. Then 
$\cro{\cdot\,}{\cdot}$ is positive
definite on
$\Pg$, negative definite on $\Kg$ and we have
\begin{equation}\label{orthog}
\cro{\Pg}{\Kg}=0.
\end{equation} 
Among others, one reason for this normalization is 
that the scalar product on
$\Pg\simeq T_{eK}(G/K)$ defined by the restriction of 
$\cro{\cdot\,}{\cdot}$ induces precisely
the
$G$-invariant Riemannian metric on
$\Hyp^n_\K=G/K$ which has pinched sectional curvature inside the 
interval $[-4,-1]$ (and constant, equal to $-1$, in the real case).

For $t\in\R$, we set $a_t=\exp(tH_0)$, so that
$$A=\{a_t,\, t\in\R\}.$$ 
We have the classical \textit{Cartan decomposition} $G=KAK$, which
actually can be slightly refined as
\begin{equation} \label{CartanG}
G=K\{a_t,\, t\ge 0\}K.
\end{equation}
When writing $g=k_1 a_t k_2$ with $t\ge 0$ according to decomposition
\eqref{CartanG}, we then have
\begin{equation}\label{disthyp}
t= \textrm{hyperbolic distance }d(gK,eK),
\end{equation} 
where $eK$ is the origin in $\Hyp^n_\K=G/K$.

\subsection{Differential forms}

In order to explain the way we shall view a differential form on a hyperbolic space,
let us proceed with some tools coming from representation theory of the groups $G$ and $K$.
First, denote as usual by
$M$ the centralizer of
$A$ in
$K$, with corresponding Lie algebra
$\Mg$, and let
$P=MAN$ be the standard minimal parabolic subgroup of $G$.  For 
$\sigma\in \widehat{M}$ and
$\lambda\in \Ag^*_\C\simeq\C$, the \textit{principal series representation}
$\pi_{\sigma,\la}$ of $G$ is
the induced representation
$$\pi_{\sigma,\la}=\Ind^G_P\left(\sigma\otimes  e^{i\lambda}\otimes 
\textbf{1}\right)$$ with
corresponding space
\begin{equation*} H^\infty_{\sigma,\tau}=\{f\in 
C^\infty(G,V_\sigma),\,  f(xma_t
n)=e^{-(i\lambda+\rho)t}
\sigma(m)^{-1}f(x),
\forall x\in G,\,\forall ma_t n\in  P \}.
\end{equation*} This $G$-action is given by left translations:
$\pi_{\sigma,\la}(g)f(x)=f(g^{-1}x)$. Moreover, if $H_{\sigma,\la}$ 
denotes the Hilbert
completion of $H^\infty_{\sigma,\la}$ with respect to the norm $\|f\|=\left\|
f|_K\right\|_{L^2(K)}$, then $\pi_{\sigma,\la}$  extends to a 
continuous representation of $G$
on $H_{\sigma,\la}$. When $\lambda\in \R$, the principal series representation
$\pi_{\sigma,\la}$ is unitary, in which case it is also irreducible, 
except maybe for
$\lambda=0$.

Next, let $(\tau, V_\tau)$ be a unitary finite dimensional 
representation of the group
$K$ (not necessarily irreducible). It is standard (\cite{Wallach73}, 
\S 5.2) that the space of
sections of the
$G$-homogeneous vector bundle
$E_\tau=G\times_K V_\tau$ can be identified with the space
$$\Gamma(G,\tau)=\{f:G\to V_ \tau,\  f(xk)=\tau(k)^{-1}f(x),\forall 
x\in G,\,\forall k\in K\}$$
of functions of (right) type $\tau$ on $G$. We define also the subspaces
$$C^\infty(G,\tau)=\Gamma(G,\tau)\cap C^\infty(G,V_\tau),
\quad\textrm{and}\quad L^2(G,\tau)=\Gamma(G,\tau)\cap L^2(G,V_\tau)$$ 
of $\Ga(G,\tau)$ which
correspond to $C^\infty$ and $L^2$ sections of $E_\tau$, 
respectively. Note that
$L^2(G,\tau)$ is the Hilbert space associated with the unitary 
induced representation
$\Ind_K^G(\tau)$ of $G$, the action being given by left translations.

For $0\le p\le dn$, let $\tau_p$ denote the $p$-th exterior product 
of the complexified
coadjoint representation ${\rm Ad}^*_\C$ of $K$ on $\Pg^*_\C$. Then 
$\tau_p$ is a unitary
representation of $K$ on $V_{\tau_p}=\wedge^p\Pg^*_\C$ and the 
corresponding homogeneous bundle
$E_{\tau_p}$ is the bundle of differential forms of degree $p$ on $G/K$.

In general, the representation $\tau_p$ is not $K$-irreducible and 
decomposes as a finite
direct sum of $K$-types:
\begin{equation}\label{dectaup}
\tau_p=\bigoplus_{\tau\in \widehat{K}} m(\tau,\tau_p)\tau,
\end{equation} where $m(\tau,\tau_p)\ge 0$ is the multiplicity of $\tau$ in
$\tau_p$ (as usual, $\widehat K$ stands for the unitary dual of the 
Lie group $K$). Let us set
$$\widehat{K}(\tau_p)=\{\tau\in \widehat{K},\, m(\tau,\tau_p)>0\},$$ 
so that \eqref{dectaup}
induces the following decomposition:
\begin{equation} \label{decEtaup}
L^2(G,\tau_p)=\bigoplus_{\tau\in \wh
K(\tau_p)}\left(L^2(G,\tau)\otimes\C^{m(\tau,\tau_p)}\right),
\end{equation} as well as its analogue when considering $C^\infty$ 
differential $p$-forms.

\subsection{The continuous part of the Plancherel formula for $L^2(G,\tau_p)$}

Let us consider an irreducible unitary representation $\tau\in 
\widehat{K}$. When restricted to
the subgroup $M$, $\tau$ is generally no more irreducible, and splits 
into a finite direct sum
$$\tau|_M=\bigoplus_{\sigma\in \widehat{M}}m( \sigma,\tau)\sigma,$$ 
where $m( \sigma,\tau)\ge
0$ is the multiplicity of $\sigma$ in $\tau|_M$ and $\widehat{M}$ stands for the unitary dual
of $M$. Let us define then
$$\widehat{M}(\tau)=\{\sigma\in \widehat{M},\, m( \sigma,\tau)>0\}.$$

The \textit{Plancherel formula} for the space $L^2(G,\tau)$ of
$L^2$ sections of the homogeneous bundle $E_\tau=G\times_K V_\tau$ 
consists in the
diagonalization of the corresponding unitary representation 
$\Ind_K^G(\tau)$ of $G$. First, we
remark that
$$L^2(G,\tau)\simeq\{L^2(G)\otimes V_\tau\}^K,$$ where the upper 
index $K$ means that we take
the subspace of $K$-invariant vectors for the right action of $K$ on 
$L^2(G)$. According to
Harish-Chandra's famous Plancherel Theorem for $L^2(G)$, the space
$L^2(G,\tau)$ splits then into the direct sum of a continuous part 
$L^2_c(G,\tau)$ and of a
discrete part $L^2_d(G,\tau)$. The latter can be expressed in terms 
of discrete series representations of $G$, but giving such a precision would be useless for our 
purpose. The former takes
the following form (see e.g.
\cite{Pedon99},
\S 3, for details):
\begin{equation} \label{planchtau}
L^2_c(G,\tau) \simeq\bigoplus_{\sigma\in
\widehat{M}(\tau)}\int^\oplus_{\Ag_+^*} d\lambda\, p_\sigma(\lambda)
H_{\sigma,\lambda}\widehat{\otimes}\Hom_K(H_{\sigma,\lambda},V_\tau) 
\end{equation} 
In this formula, $d\lambda$ is the Lebesgue measure on 
$\Ag^*_+\simeq
(0,+\infty)$, $p_\sigma(\lambda)$ is the Plancherel density 
associated with $\sigma$ and
$\Hom_K(H_{\sigma,\lambda},V_\tau)$ is the vector space of
$K$-intertwining operators from $H_{\sigma,\lambda}$ to $V_\tau$, on 
which $G$ acts trivially.
This space is non zero (since $\si\in\wh M(\tau)$) but finite 
dimensional (since every
irreducible unitary representation of $G$ is admissible).

By combining formulas \eqref{decEtaup} and \eqref{planchtau} we get 
the  following result.

\begin{prop} \label{planchtaup}
The continuous part of the Plancherel formula for 
$L^2(G,\tau_p)$ is given by:
\begin{equation*} L^2_c(G,\tau_p) \simeq
\bigoplus_{\tau\in\widehat{K}(\tau_p)}
\left(\bigoplus_{\sigma\in\widehat{M}(\tau)}
\int^\oplus_{\Ag_+^*} d\lambda\, p_\sigma(\lambda)
H_{\sigma,\lambda}\widehat{\otimes}\Hom_K(H_{\sigma,\lambda},V_\tau)
\right)\otimes\C^{m(\tau,\tau_p)}.
\end{equation*}
\end{prop}

\subsection{The spectrum of the Hodge-de~Rham Laplacian}
\label{specHDR}

The \textit{Hodge-de~Rham Laplacian} $\De_p=dd^*+d^*d$ acts on $C^\infty$ differential
$p$-forms on $\Hyp^n_\K=G/K$, i.e. on members of the space
$C^\infty(G,\tau_p)$. Actually, this operator is realized by the 
action of the Casimir element
$\Om_\Gg$ of the universal enveloping algebra $U(\Gg)$ of $\Gg$. More 
precisely, keeping
notation \eqref{bilin}, let
$\{Z_i\}$ be any basis for
$\Gg$ and $\{Z^i\}$ the corresponding basis of $\Gg$ such that 
$\cro{Z_i}{Z^j}=\de_{ij}$. The
Casimir operator can be written as
\begin{equation}\label{Casimir}
\Om_\Gg=\sum_i Z_i Z^j.
\end{equation} We can view $\Omega_\Gg$ as a $G$-invariant 
differential operator acting on
$C^\infty(G,\tau_p)$, and have then the well-known identification 
(\textit{Kuga's formula}, see
\cite{BW00}, Theorem~II.2.5)
\begin{equation}\label{Kuga}
\Delta_p=-\Omega_\Gg.
\end{equation}

We shall denote also by $\De_p$ the unique self-adjoint extension of the Hodge-de~Rham operator
from compactly supported smooth differential forms to $L^2$ differential forms on
$\Hyp^n_\K=G/K$. Let us recall that the nature of its spectrum is well known:

\begin{thm} \label{specGK}
\begin{enumerate}
\item If $p\not= \frac{dn}{2}$, the $L^2$ spectrum of $\De_p$ is absolutely continuous, of the form
$[\al_p,+\infty)$ with $\al_p\ge 0$.
\item If $p=\frac{dn}{2}$ (with $dn$ even), one must add the sole discrete eigenvalue
$0$, which occurs with infinite multiplicity.
\item We have $\al_p=0$ if and only if $\K=\R$ and $p=\frac{n\pm1}{2}$. In particular, the discrete
eigenvalue $0$ occuring in middle dimension $p=\frac{dn}{2}$ is always spectrally isolated.
\end{enumerate}
\end{thm}

In this result, assertion (1) is essentially \pref{planchtaup}, assertion (2) is true for
any general $G/K$ and can be found e.g. in \cite{Borel85}, \cite{Pedon97} or \cite{Olbrich}, and
assertion (3) follows from results in \cite{BW00}
and \cite{VZ84}, as noticed by J.~Lott in \cite{Lott}, \S~VII.B (see also our \tref{valalphap}).

Moreover, the exact value of $\al_p$ can be calculated with the help of some more
representation theory. Let us elaborate. Thanks to \eqref{Kuga}, in order to investigate the
continuous $L^2$ spectrum of
$\Delta_p$ and thus to compute $\al_p$,  it is enough to
consider the action of the Casimir operator $\Omega_\Gg$ on the 
right-hand side of the
Plancherel formula given in \pref{planchtaup} and, specifically, on 
each elementary
component 
$H_{\sigma,\lambda}\widehat{\otimes}\Hom_K(H_{\sigma,\lambda},V_\tau)$.

The action of $\Omega_\Gg$ on $\Hom_K(H_{\sigma,\lambda},V_\tau)$ 
being trivial, the problem
reduces to study its effect on $H_{\si,\la}$, and even on 
$H^\infty_{\si,\la}$, by density. But
since
$\Om_\Gg$ is a central element in the enveloping algebra of $\Gg$, it 
acts on the irreducible
admissible representation
$H^\infty_{\si,\la}$ by a scalar $\om_{\si,\la}$. More precisely, 
let $\mu_\sigma$ be the
highest weight of $\sigma \in \widehat{M}$ and $\delta_\Mg$ be the 
half sum of the positive
roots of $\Mg_\C$ with respect to a given Cartan subalgebra. Then $\si(\Om_\Mg)=-c(\si)\Id$,
where the Casimir value of $\si$ is given by
\begin{equation} \label{c(si)}
c(\si)=\cro{\mu_\si}{\mu_\sigma+2\delta_\Mg}\ge 0.
\end{equation} 
Using for instance
\cite{Knapp86}, Proposition~8.22 and Lemma~12.28, one easily checks that
\begin{equation}\label{CasSP}
\Omega_\Gg=\omega_{\sigma,\lambda} \Id
\textrm{ on }H^\infty_{\sigma,\lambda},
\end{equation}
where
$$\omega_{\sigma,\lambda}=-\left(\lambda^2+\rho^2-c(\si)\right).$$
Thus \eqref{CasSP},
\eqref{Kuga} and 
\pref{planchtaup} show that the action
of $\De_p$ on (smooth vectors of) $L^2(G,\tau_p)$ is diagonal, a fact 
which allows us to
calculate the continuous $L^2$ spectrum of $\De_p$.

In order to state this, set
$$\wh{M}(\tau_p)=\bigcup_{\tau\in\wh{K}(\tau_p)}\wh{M}(\tau),$$ and 
denote by $\si_{\max}$ one
of the (possibly many) elements of $\wh{M}(\tau_p)$ such that
$c(\si_{\max})\ge c(\si)$ for any $\si\in\wh{M}(\tau_p)$. Our 
discussion implies immediately
the following result.

\begin{prop} \label{bottom}
The continuous $L^2$ spectrum of the Hodge-de~Rham 
Laplacian $\Delta_p$ is $[\alpha_p,+\infty)$, where 
\begin{equation}\label{alphap}
\alpha_p=\rho^2-c(\si_{\max}).
\end{equation}
\end{prop}

With a case-by-case calculation, the previous formula gives the explicit value of $\al_p$ (at
least in theory; in the case $\K=\H$, identifying the representations $\si_{\max}$ is quite
awkward, see \cite{Pedon}). For instance,
$\al_0$ equals
$\rho^2$ for any
$\Hyp^n_\K$, since
$\si_{\max}$ must be the trivial representation (this well-known fact can be proved also by other
arguments). For general
$p$, we collect the known results in the following theorem. Observe that we can restrict to $p\le
dn/2$, since $\al_{dn-p}=\al_p$ by Hodge duality.

\begin{thm} \label{valalphap}
Let $p\le \frac{dn}{2}$.
\begin{enumerate}
\item If $\K=\R$ (see \cite{Donnelly81},
\cite{Pedon98a}), then $\al_p=\left(\frac{n-1}{2}-p\right)^2$. 
\item If $\K=\C$ (see \cite{Pedon99}), then 
\begin{equation*}
\al_p=
\begin{cases}
(n-p)^2&\textrm{if }p\not= n,\\
1&\textrm{if }p=n.
\end{cases} 
\end{equation*}
\item If $\K=\H$ (see \cite{Pedon}), then
\begin{equation*}
\al_p=
\begin{cases} 
(2n+1)^2&\text{if }p=0,\\ 
(2n-p)^2+8(n-p)&\text{if }1\le p\le[\frac{4n-1}{6}],\\
(2n+1-p)^2&\text{if }[\frac{4n-1}{6}]+1\le p\le n,\\ 
(2n-p)^2&\text{if }n+1\le p\le 2n-1,\\ 
1&\text{if }p=2n.
\end{cases}
\end{equation*} 
\end{enumerate}
\end{thm}

To our knowledge, the value of $\al_p$ in the exceptional case
$\Hyp^2_\O$ is still unknown.

\subsection{The action of the Hodge-de~Rham Laplacian on 
$\tau_p$-radial functions}

For any finite dimension representation $\tau$ of $K$, let us introduce the space of 
\textit{smooth $\tau$-radial functions on $G$}:

\begin{equation}\begin{split} C^\infty(G,\tau,\tau)&=\{F\in C^\infty(G,\End V_{\tau}),\\
&F(k_1gk_2)=\tau_p(k_2)^{-1} F(g)\tau_p(k_1)^{-1},
\forall g\in G,\, \forall k_1,k_2\in K \}.
\end{split}\end{equation}

Our aim is to calculate the action of the Laplacian $\De_p$ on $C^\infty(G,\tau_p,\tau_p)$ (the
reason will be given in next subsection). Because of the Cartan decomposition 
\eqref{CartanG}, it is clear
that any
$\tau$-radial function on
$G$ is entirely determined by its restriction to the semigroup 
$\{a_t,\,t\ge 0\}$. Hence it is
sufficient to calculate the value of $\De_p F(a_t)$ for any $F\in 
C^\infty(G,\tau_p,\tau_p)$ and
any $t\ge 0$.

Because of \eqref{Casimir} we have 
$$\Omega_\Gg=\Omega_\Pg-\Omega_\Kg=\sum_i X_i^2-\sum_i
Y_i^2$$ if we choose bases $\{X_i\}$ of $\Pg$ and $\{Y_i\}$ of $\Kg$ 
which satisfy
respectively $\langle X_i,X_j\rangle=\delta_{ij}$ and $\langle 
Y_i,Y_j\rangle=-\delta_{ij}$.
On the spaces $C^\infty(G,\tau_p)$ and $C^\infty(G,\tau_p,\tau_p)$, we thus get
\begin{equation}\label{Weitz}
\Delta_p=-\Omega_\Gg=-\Omega_\Pg+\tau_p(\Omega_\Kg),
\end{equation} where $\tau_p(\Omega_\Kg)$ is a zero order 
differential operator  which is
diagonal, since $\tau(\Omega_\Kg)$ is scalar for each $\tau\in 
\widehat K$, namely
\begin{equation*}
\tau(\Om_\Kg)=-c(\tau)\Id=-(\mu_\tau\vert\mu_\tau+2\de_\Kg)\Id,
\end{equation*} with notations that are analogous to the ones used in 
\eqref{c(si)}. Notice that
\eqref{Weitz} is exactly the well-known Bochner-Weitzenb\"ock 
formula (see \eqref{BochnerW}), since $-\Omega_\Pg$
coincides with the Bochner Laplacian $\nabla^*\nabla$ (see e.g. 
\cite{BOS94}, Proposition~3.1).

Now, reminding \eqref{defn} and \eqref{orthog}, let $\Lg_1$ and 
$\Lg_2$ be the orthogonal
projections of the root subspaces $\Gg_\al$ and $\Gg_{2\al}$ on
$\Kg$ with respect to the Cartan decomposition \eqref{Cartang}, so 
that we have the orthogonal
splitting
\begin{equation}\label{deck}
\Kg=\Mg\oplus \Lg_1\oplus\Lg_2.
\end{equation} (Remark that $\Lg_2$ reduces to zero if $\K=\R$.) Let
$\{Y_{1,r}\}_{r=1}^{d(n-1)}$ and $\{Y_{2,s}\}_{s=1}^{d-1}$ denote the 
subsystems of the basis
$\{Y_i\}$ of $\Kg$ which are bases for $\Lg_1$ and $\Lg_2$, respectively.

We have then the following result.

\begin{prop} \label{Bochner}
If $F\in C^\infty(G,\tau_p,\tau_p)$, then for any $t\ge 0$ we have
\begin{equation*}
\De_p F(a_t)=-\Om_\Pg F(a_t)+\tau_p(\Om_\Kg)F(a_t),
\end{equation*} where
\begin{multline}\label{eqBochner}
\Omega_\Pg F(a_t)=
\frac{d^2}{dt^2}F(a_t)+[d(n-1)\coth t+2(d-1)\coth 
2t]\frac{d}{dt}F(a_t)\\ +(\coth
t)^2\sum_{r=1}^{d(n-1)} \tau_p(Y_{1,r}^2)F(a_t)+(\sinh 
t)^{-2}F(a_t)\sum_{r=1}^{d(n-1)}
\tau_p(Y_{1,r}^2)\\ -2(\sinh t)^{-1}(\coth t)\sum_{r=1}^{d(n-1)}
\tau_p(Y_{1,r})F(a_t)\tau_p(Y_{1,r})\\ +(\coth 2t)^2\sum_{s=1}^{d-1}
\tau_p(Y_{2,s}^2)F(a_t)+(\sinh  2t)^{-2}F(a_t)\sum_{s=1}^{d-1} 
\tau_p(Y_{2,s}^2)\\ -2(\sinh
2t)^{-1}(\coth 2t)\sum_{s=1}^{d-1} \tau_p(Y_{2,s})F(a_t)\tau_p(Y_{2,s}).\\
\end{multline}
\end{prop}

\begin{proof} It remains only to show formula \eqref{eqBochner}, whose proof 
is standard and can be found e.g. in
\cite{Wallach73}, \S 8.12.6.
\end{proof}
 
\subsection{The resolvent of the Hodge-de~Rham Laplacian and the associated Green kernel}
\label{sresolv}

It is well-known that all kernels $K(x,y)$ of functions of the (positive)
Laplace-Beltrami operator $\De_0$ on a symmetric space $G/K$ only depend
on the Riemannian distance: $K(x,y)=k(d(x,y))$. In other words, because of \eqref{disthyp}
and the $G$-invariance of the distance, they can be considered as radial (i.e. 
bi-$K$-invariant) functions on $G$.
In the case of our bundle of differential forms, kernels of operators related
to the Hodge-de~Rham Laplacian $\De_p$ will naturally be $\tau_p$-radial 
functions on $G$ (see e.g.
\cite{CM82}). In particular, for $s\in\C$ with
$\re s>0$, consider the \textit{Green kernel} $G_p(s,\cdot)$ of the \textit{resolvent}
$$R_p(s)=(\De_p-\al_p+s^2)^{-1}.$$
By definition, it solves the differential equation
\begin{equation}\label{eqGreen}
(\De_p-\al_p+s^2)G_p(s,\cdot)=\de_e.
\end{equation}
Therefore, when $\re s>0$ the Green kernel $G_p(s,\cdot)$ is a Schwartz $\tau_p$-radial 
function on $G^0=G\smallsetminus\{e\}$, i.e.
a member of the space
\begin{equation*}
\begin{split} {\mathcal S}(G^0,\tau_p,\tau_p)=\{F\in & 
\,C^\infty(G^0,\tau_p,\tau_p)\,:\,
\forall D_1,D_2\in U(\Gg),\ \forall N\in\N,\\ &\sup_{t>0}\|F(D_1:a_t:D_2)\|_{\End
V_{\tau_p}}(1+t)^N e^{\rho t}<+\infty\},
\end{split}
\end{equation*} 
where we use the classical Harish-Chandra notation $F(D_1:a_t:D_2)$ for the two
sided derivation of $F$ at $a_t$ with respect to the elements $D_1$ and $D_2$ of the universal
enveloping algebra $U(\Gg)$. 

For convenience, we shall often use the alternative notation
\begin{equation}\label{notaGreen}
g_p(s,t)=G_p(s,a_t),
\end{equation}
defined for $\re s>0$ and $t\ge 0$ (with a singularity at $t=0$).

\subsection{Hyperbolic manifolds}

Throughout this paper, $\Ga$ will denote any torsion-free discrete subgroup of $G$, so that
the quotient $\Ga\backslash G/K$ is a \textit{hyperbolic manifold}, 
i.e. a complete Riemannian locally
symmetric space with strictly negative curvature. 
We define
$\de(\Ga)$ to be the \textit{critical exponent} of the Poincar\'e series associated with $\Ga$,
i.e. the nonnegative number
$$\delta(\Gamma)=\inf\{s\in\R \textrm{\ such\ that\ } \sum_{\gamma\in\Gamma} e^{-s
d(x,\gamma y)} <+\infty\},$$
where $(x,y)$ is any pair of points in $\Hyp^n_\K$ (for instance, $x=y=eK$) and $d$ is the
hyperbolic distance. It is easy to check that 
$$0\le\de(\Ga)\le 2\rho=h,$$
and it is also known that equality holds if $\Ga$ has finite covolume (actually the converse is
true when $\K=\H$ or $\O$, see \cite{corlette}, Theorem~4.4). 
This critical exponent has been extensively studied. For instance, one of the most striking
results says that when
$\Gamma$ is geometrically finite, then $\delta(\Gamma)$ is the Hausdorff dimension of the limit
set $\Lambda(\Gamma)\subset S^\infty$, where the sphere at infinity $S^\infty=\partial
\Hyp^n_\K\simeq K/M\simeq{\bS}^{dn-1}$ is endowed with its natural Carnot structure (see the 
works of S.~Patterson, D.~Sullivan, C.~Yue (\cite{patterson}, \cite{sullivand},
\cite{yue}). 
Regarding the \textit{limit set} $\Lambda(\Ga)$, let us recall that it is defined as the
set of accumulation points of any $\Ga$-orbit in the natural compactification
$\Hyp^n_\K\cup S^\infty$.

In some cases, we also consider more general locally symmetric spaces $\Ga\backslash G/K$ of the
noncompact type (i.e. with nonpositive sectional curvature, of any rank). All
definitions and results 
described above extend to that situation (except the relationship between $\de(\Ga)$ and
$\Lambda(\Ga)$, which is still being investigated, see \cite{albu}, \cite{QuintG}, \cite{Quint}).

Finally, if $X$ stands for any complete Riemannian manifold, we let 
$\lambda_0^p(X)$ be the bottom of the $L^2$ spectrum of $\De_p$ on 
$X$. In other words,
$$\lambda_0^p(X)=\inf_{u\in C^\infty_0(\wedge^p T^*X)}
\frac{(\De_p u|u)_{L^2}} {\|u\|_{L^2}^2}.$$
With notation of \tref{specGK}, we thus have $\la_0^p(\Hyp^n_\K)=\al_p$ when
$p\not=\frac{dn}{2}$. Note also that, when $p=0$, this definition reduces to
$$\lambda_0^0(X)=\inf_{u\in C^\infty_0(X)}
\frac{\|du\|_{L^2}^2} {\|u\|_{L^2}^2}.$$

\subsection{Some remarks about the generalization of \tref{SC} to other Riemannian manifolds}

Hyperbolic spaces admit natural generalizations. Namely, they can be viewed both as a particular
class of symmetric spaces of noncompact type and as a particular class of \textit{harmonic $AN$
groups} (also called \textit{Damek-Ricci spaces}). The latter are Einstein manifolds which are not
symmetric (except for the hyperbolic spaces) but their analysis is quite similar to the one of
hyperbolic spaces (see \cite{ADY}). 

For these two families of manifolds, the proof of \tref{SC} can be adapted to get information on
the bottom of the spectrum of the Laplacian $\De_0$ defined on some quotient by a
discrete torsion-free subgroup
$\Ga$. Indeed, in both cases one has at his disposal the key ingredient, that is, estimates for
the Green kernel (see Theorem~4.2.2 in
\cite{AJ} and Theorem~5.9 in \cite{ADY}, respectively). 

In the case of Damek-Ricci spaces $AN$,
the result reads exactly as in \tref{SC}, provided we replace
$2\rho$ by the homogeneous dimension of $N$ (in both cases, these numbers represent the
exponential rate $h$ of the volume growth). 

As concerns locally symmetric spaces $\Ga\backslash G/K$, the statement is not as sharp as in
\tref{SC}, since it provides in general only bounds for $\la_0^0(\Ga\backslash G/K)$. Let us
elaborate.

Take the Lie groups $G$ and $K$ as in Section~\ref{notation},
except that $G/K$ can be now of any rank $\ell\ge 1$, which means that $\Ag\simeq\R^\ell$.
Let us introduce some more notation. First, we
have an inner product on all $\Gg$ by modifying the symmetric bilinear form \eqref{bilin} as
follows: 
\begin{equation}\label{bilin2}
\cro{X}{Y}=-B(X,\theta Y)\qquad\forall X,Y\in\Gg,
\end{equation}
and we denote by $\norm{\cdot}$ the corresponding norm. The restriction of
\eqref{bilin2} to $\Pg$ induces a $G$-invariant Riemannian metric on
$G/K$ of (non strictly if $\ell>1$) negative curvature. 

For any element $x\in G$, define $H(x)$ to be the unique element in
the closure $\overline{\Ag_+}$ of the positive Weyl chamber in
$\Ag$ so that
$$x=k_1\exp H(x) k_2$$
reflects the Cartan decomposition of $x$ (the analogue of \eqref{CartanG}). The half-sum
$\rho\in\Ag^*$ of positive roots of the pair $(\Gg,\Ag)$ cannot be considered as a real number
anymore. Nevertheless, we can still view it as a member of $\Ag$ via \eqref{bilin2}, and it
should be noted that
$$\la_0^0(G/K)=\norm{\rho}^2,\qquad h=2\norm{\rho},$$
as well as
$$0\le\de(\Ga)\le 2\norm{\rho}.$$

Using Theorem~4.2.2 in \cite{AJ}, E.~Leuzinger has obtained the following result (see
\cite{Leuzinger}).

\begin{thm} Let $G/K$ be any noncompact Riemannian symmetric space, let $\Ga$ be a discrete
torsion-free subgroup of $G$, and set $\rho_{\min}=\inf_{H\in \overline{\Ag_+}} \langle
\rho,H\rangle /\norm{H}\,$ (so that $\rho_{\min}\le\norm{\rho}$, with equality in the rank one
case). 
\begin{enumerate}
\item If $\delta (\Gamma)\le \rho_{\min}$, then
$\lambda_0^0(\Gamma\backslash G/K)=\norm{\rho}^2$.
\item If $\delta (\Gamma)\in[ \rho_{\min},\norm{\rho}]$, then
$$\norm{\rho}^2-(\delta (\Gamma)-\rho_{\min})^2
\le\lambda_0^0(\Gamma\backslash G/K)\le \norm{\rho}^2.$$
\item If $\delta (\Gamma)\ge \norm{\rho}$, then
$$\max\{\norm{\rho}^2-(\delta (\Gamma)-\rho_{\min})^2,0\}
\le\lambda_0^0(\Gamma\backslash G/K)\le \delta
(\Gamma)(2\norm{\rho}-\delta (\Gamma)).$$
\end{enumerate}
\end{thm}

Actually, we have a better expression in terms of a modified critical
exponent. The proof is underlying in \cite{Leuzinger}.

\begin{thm} Let $G/K$ be any noncompact Riemannian symmetric space, and let $\Ga$ be a discrete
torsion-free subgroup of $G$.
Define $\tilde\delta (\Gamma)$ to be the critical exponent of
the Poincar\'e series $$\sum_{\gamma\in\Gamma} e^{-sd(x,\gamma y)
-\rho(H(\gamma))}.$$ 
(This definition does not depend on the points $x,y\in G/K$.)
Then $$\lambda_0^0(\Gamma\backslash G/K)=\norm{\rho}^2-\tilde\delta (\Gamma)^2.$$
\end{thm}

In addition to this statement, one other reason to introduce our modified critical exponent is
motivated by the following observation. Suppose that $\Gamma$ is Zariski dense in
$G$. 
In \cite{Quint}, J.-F.~Quint defines a function $\Phi_\Gamma\,:\,
\overline{\Ag_+}\to \R\cup\{-\infty\}$ which measures the growth
of $\Gamma$ in the direction of $H\in \overline{\Ag_+}$, and he shows
that this function is concave. According to Corollary~5.5 in
\cite{Quint}, we have a link between the growth indicator $\Phi_\Ga$ and our modified
critical exponent $\tilde\de(\Ga)$, namely:
$$ \tilde\delta (\Gamma)=\inf_{\substack{H\in \overline{\Ag_+}\\ \norm{H}=1}} \left(\langle
\rho,H\rangle+\Phi_\Gamma(H)\right).$$

\section{The resolvent associated with the Hodge-de~Rham Laplacian on hyperbolic spaces}

In our investigation of the bottom of the differential $p$-form spectrum on a hyperbolic manifold
$\Ga\backslash\Hyp^n_\K$,
the key step consists in the careful analysis of the resolvent 
$R_p(s)=(\De_p-\al_p+s^2)^{-1}$ associated with the covering space $\Hyp^n_\K$, and this goes
through an estimate of the corresponding Green kernel $G_p(s,\cdot)$.

As a matter of fact, these estimates will partly be obtained by comparing $G_p(s,\cdot)$ to the
scalar Green kernel $G_0(s,\cdot)$. We thus begin with the following result, whose proof is
standard but will be recalled here, since we need to emphasize some of its ingredients.
We retain notation from previous sections and particularly from \sref{sresolv}.

\begin{prop} \label{prop1}
For any $(\Ga x,\Ga y)\in(\Ga\backslash\Hyp^n_\K)\times (\Ga\backslash\Hyp^n_\K)$ 
and any $s\in\C$ with $\re s>0$, let 
$$g^*_0(s,\Ga x,\Ga y)=\sum_{\gamma\in \Gamma}g_0(s,d(x,\gamma y))$$ 
be the pull-back of the
Green kernel from $\Ga\backslash\Hyp^n_\K$ to $\Hyp^n_\K$.
Then, for any $s>0$, $g^*_0(s,\Ga x,\Ga y)$ behaves as the Poincar\'e series
$$\sum_{\ga\in\Ga}e^{-(s+\rho)d(x,\ga y)}.$$
\end{prop}

\begin{proof}
We first observe that the Green kernel $g_0(s,\cdot)$ defined by \eqref{notaGreen} (and
corresponding to the resolvent $R_0(s)=(\De_0-\rho^2+s^2)^{-1}$)
can be explicitly expressed as a hypergeometric function (see for
instance \cite{faraut},
\cite{Miawa}, \cite{ADY}). Indeed, by \eqref{eqGreen} and \pref{Bochner}, it must solve the
Jacobi type differential equation
$$g_0''(s,r)+[(dn-1)\coth r+(d-1)\tanh r] g_0'(s,r)+(\rho^2-s^2) g_0(s,r)=0,$$
where differentiation is meant with respect to the second variable $r$.
Letting 
$$u(s,-(\sinh r)^2)=g_0(s,r),$$ 
we see that the function $u$
solves the hypergeometric equation
$$ x(1-x)
u''(s,x)+\left(\frac{dn}{2}-\frac{d(n+1)}{2}x\right)u'(s,x)-\frac{\rho^2-s^2}{4}
u(s,x)=0.$$
Since the resolvent
$R_0(s)$ acts continuously on $L^2(\Hyp_\K^n)$
and since we must have the following standard behaviour:
\begin{equation}\label{eqlowerbound0}
g_0(s,r)\underset{r\to 0}{\simeq}
\left\{
\begin{array}{ll}
\frac{r^{2-dn}}{\vol(\bS^{dn-1})} & {\rm
if\ } dn>2, \\
-\frac{1}{2\pi}\log r &  {\rm if\ } dn=2, \\
\end{array}\right.
\end{equation}
by using
theorem 2.3.2 in \cite{AAR} we find that
$$u(s,x)=f_{n,d}(s)\, (2x)^{-(s+\rho)/2}
{\,}_2F_1\left(\frac{s+\rho}{2},\frac{s+1}{2}-\frac{d(n-1)}{4},s+1,x^{
-1}\right),$$
where ${}_2F_1$ is the classical Gauss hypergeometric function and
$$
f_{n,d}(s)=2^{d-2} \pi^{-(dn-1)/2} \frac
{\Gamma\left(\frac{s+\rho}{2}\right)\Gamma\left(
s+\frac{d(n-1)}{2}\right)} {\Gamma\left( s+1
\right)\Gamma\left(\frac{s}{2}+\frac{d(n-1)}{4}\right)}.$$

>From these explicit formulas, we deduce important facts. Firstly,
the resolvent
$$R_0(s)\,:\,C_0^\infty(\Hyp_\K^n)\longrightarrow C^\infty(\Hyp_\K^n),$$
which is a priori defined for $\re s>0$, has a meromorphic extension
to the complex plane and has a
holomorphic extension to the half-plane $\re s>-\frac{d(n-1)}{2}$.

Secondly, we can estimate the function
$g_0(s,r)$ for large values of $r$ (see also \cite{LR}). On the one hand, for every $s\in \C$ such
that $\re s>-\frac{d(n-1)}{2}$, there is a positive constant $c_1(s)$ such that
\begin{equation}\label{estimgreenmero}
\forall r\geq 1, \ \ |g_0(s,r)|\leq c_1(s)\, e^{-(\re s+\rho) r}.
\end{equation}
On the other hand, when $s$ is a positive real number, $g_0$ is a positive
real function and it can be bounded from below: there exists a positive
constant $c_2(s)$ such that
\begin{equation}\label{eestimgreen}
\forall r\geq1,\ \ c_2(s)\,e^{-(s+\rho)r}\leq g_0(s,r).
\end{equation}
The result immediately follows.
\end{proof}

Assume that $\delta(\Gamma)< \rho$.
As was noticed by Y. Colin de Verdi\`ere in \cite{CdV}, the estimate
\eqref{estimgreenmero} implies also that
the resolvent of the Laplacian $\De_0$ on $\Gamma \backslash\Hyp_\K^n$ 
has a holomorphic continuation to the
half-plane 
$$\{s\in\C\,:\,\re s>  \min( -\tfrac{d(n-1)}{2},\delta(\Gamma)-\rho) \}.$$ 
According to a well known principle
of spectral theory (\cite{RS}, Theorem XIII.20), 
we thus get:

\begin{cor}  \label{remcdv}
If $\de(\Ga)<\rho$, the $L^2$ spectrum of
the Laplacian $\De_0$ on $\Gamma \backslash\Hyp_\K^n$ is absolutely continuous.
\end{cor}

Now we turn to the general case of differential forms.
Reminding  notation from \sref{specHDR}, we define
$\Sigma$ to be the (minimal) branched cover of $\C$ such that the functions
$s\mapsto\sqrt{s^2-c(\sigma)+c(\si_{\max})}$ are holomorphic on
$\Sigma$ for all $\sigma\in \widehat{M}(\tau_p)$. This cover is realized as follows: let
$\si_1,\ldots,\si_r$ denote the distinguished representatives of the $\si$'s in
$\widehat{M}(\tau_p)$ such that
$c(\si)\not=c(\si_{\max})$. Then
\begin{equation}\label{Sigma}
\Sigma=\{\hat s=(s,y_1,\ldots,y_r)\in \C^{r+1}\,:\,y_i^2=s^2+c(\sigma_{max})-c(\sigma_i),
\ \forall i=1,\ldots r\}.
\end{equation}
$\Sigma$ contains naturally a copy of the half-plane
$$\C_+=\{s\in \C,\, \re s>0\},$$
namely
\begin{equation}\label{Cplus}
\C_+\equiv\{\hat s=(s,y_1,\ldots,y_r)\in \Sigma\,:\,\re s>0,\ \re y_i>0,\ \forall i=1,\ldots r\},
\end{equation}
and we let
$\overline{\C_+}$ stands for its closure in $\Sigma$.
Also, we shall still denote by $s\,:\,\Sigma\rightarrow \C$ the holomorphic extension
of the function $s$ from $\C_+$ to $\Sigma$.
Finally, if $\hat s=(s,y_1,\cdots,y_r)\in\Sigma$ we set
$$h(\hat s)=\min\{\re s,\re y_1,\ldots,\re y_r\},$$
and we recall that we have put $G^0=G\priv\{e\}$.

\begin{prop} There exists a function $F_p(\hat s,x)$ defined on $\Si\times G^0$ such that:
\begin{enumerate}
\item the map $x\mapsto F_p(\hat s,x)$ belongs to $C^\infty(G^0,\tau_p,\tau_p)$;
\item $\hat s\mapsto F_p(\hat s,x)$ is meromorphic on $\Si$ and holomorphic on 
$\Si\priv\cN$, where
$\cN$ is a discrete subset of $\Si\smallsetminus\overline{\C_+}$;
\item $(\De_p-\al_p+s^2)F_p(\hat s,\cdot)=0$;
\item for any $\hat s\in\Si\smallsetminus\cN$, there is a constant 
$A(\hat s)>0$ such that
\begin{equation}\label{upperpSigma}
\forall t>1,\ \norm{F_p(\hat s,a_t)}_{\End V_{\tau_p}}\le A(\hat s)  e^{-\left(\rho +h(\hat s)
\right)t};
\end{equation}
\item for any $s\in\C_+$, there is a constant 
$A(s)>0$ such that
\begin{equation}\label{upperp}
\forall t>1,\ \norm{F_p(s,a_t)}_{\End V_{\tau_p}}\le A(s)  e^{-\left(\rho +\re s\right)t}.
\end{equation}
\end{enumerate}
\end{prop}

\begin{proof} For $s\in \C_+$ and $t>0$, define
$$v_p(s,t)=(\sinh t)^{d(n-1)/2}(\sinh 2t)^{(d-1)/2} G_p(s,a_t).$$ 
Using \eqref{eqGreen},
\pref{Bochner} and the standard behaviour of hyperbolic functions, we see that $v_p$ must solve
the differential equation
\begin{equation}\label{Lp}
\left(-\frac{d^2}{dt^2}+\rho^2-\alpha_p+s^2+D+W(e^{-t})\right)v_p(s,t)=0
\end{equation} on $\R^*_+$, where $W\,:\, \{z\in\C, |z|<1\}\to \End V_{\tau_p}$ is a holomorphic 
function vanishing at
$0$:
$$W(z)=\sum_{l=1}^\infty w_l z^l,\quad \text{with } w_l\in \End V_{\tau_p},$$ 
and
\begin{eqnarray*} D&=&-\sum_{r=1}^{d(n-1)} \tau_p(Y_{1,r}^2)-\sum_{s=1}^{d-1}
\tau_p(Y_{2,s}^2)+\tau_p(\Omega_\Kg)\\ &=&-\tau_p(\Omega_\Lg)+\tau_p(\Omega_\Kg)\\
&=&\tau_p(\Omega_\Mg)\qquad \qquad (\text{since } \Kg=\Mg\oplus\Lg)\\ &=&\bigoplus_{\sigma\in
\widehat M(\tau_p)}
\bigoplus_{l=1}^{m(\sigma,\tau_p)} \si(\Om_\Mg)\\ &=&\bigoplus_{\sigma\in \widehat M(\tau_p)}
\bigoplus_{l=1}^{m(\sigma,\tau_p)} [-c(\sigma)\Id_{V_\si}].
\end{eqnarray*} For convenience, let $L_p+s^2$ be the differential operator defined  by the
parentheses in the left hand-side of \eqref{Lp}. Recall from \eqref{bottom} that we have
$\alpha_p=\rho^2-c(\si_{\max})$. Hence, if we put
\begin{equation}\label{defE} E=\bigoplus_{\sigma\in \widehat M(\tau_p)}
\bigoplus_{l=1}^{m(\sigma,\tau_p)} [c(\si_{\max})-c(\sigma)]\Id_{V_\si},
\end{equation} we can rewrite $L_p$ as
$$L_p=-\frac{d^2}{dt^2}+E+W(e^{-t}).$$ 
By definition of $\Si$, the function $\hat
s\mapsto
\sqrt{E+s^2}$ is holomorphic on $\Si$. Thus, the equation $(L_p+s^2)v=0$ is of Fuchsian type, and
we can look for a  solution of the form
$$v_p(\hat s,t)=e^{-t \sqrt{E+s^2}}\sum_{l=0}^\infty a_l(\hat s)e^{-l t},$$ 
with coefficients
$a_l(\hat s)$ recursively defined by the formulas
\begin{align*} a_0(\hat s)&=\Id_{V_{\tau_p}},\\ l\left[2 
\sqrt{E+s^2}+l\,\Id_{V_{\tau_p}}\right]a_l(\hat s)&=\sum_{k=1}^l w_k\,a_{k-l}( 
\hat s).
\end{align*} 
Denote by $\cN$ the set consisting of the $\hat s\in 
\Sigma$ such that $2 \sqrt{E+s^2}+l\,{\rm Id}_{V_{\tau_p}}$ is a non invertible operator for some
$l\in\N^*$. Then $\cN$ is a discrete subset of $\Si\smallsetminus\overline{\C_+}$ and we obtain a
meromorphic map
$$\hat s\mapsto v_p(\hat s,\cdot)\in C^\infty((0,+\infty),\End V_{\tau_p})$$ 
which satisfies the
following properties:
\begin{itemize}
\item $\hat s\mapsto v_p(\hat s,\cdot)$ is holomorphic on $\Si\smallsetminus\cN$;
\item $v_p$ solves the differential equation 
$(L_p+s^2)v_p(\hat s,\cdot)=0$ on $\R^*_+$;
\item $v_p(\hat s,t)=e^{-t \sqrt{E+s^2}}[ {\rm  Id}_{V_{\tau_p}}+\text{O}(e^{-t})]$ as
$t\to+\infty$.
\end{itemize} Finally, letting $F_p$ be defined on $\Si\times G^0$ by
$$F_p(\hat s,a_t)=(\sinh t)^{-d(n-1)/2}(\sinh 2t)^{-(d-1)/2} v_p(\hat s,a_t),$$ 
and reminding formula
\eqref{rho}, we get the statements of our proposition. In particular, remark that 
\eqref{upperpSigma} follows from the estimate
\begin{equation}\label{estinterm} F_p(\hat s,a_t)\underset{t\to+\infty}{=}e^{-t
(\sqrt{E+s^2}+\rho)}\left[ {\rm Id}_{V_{\tau_p}}+\text{O}(e^{-t})\right],
\end{equation}
and that we deduce \eqref{upperp} by observing that, if $\hat s=(s,y_1,\ldots,y_r)\in\Sigma$, 
we have $\re y_j>\re s$ for all
$j$ on $\C_+$. In other words, 
\begin{equation}\label{h}
h(\hat s)=\re s \quad\textrm{on}\ \C_+.
\end{equation}
\end{proof}

Actually, the function $F_p$ we introduced in the proposition is in  some sense a multiple (in the
variable
$\hat s$) of the Green kernel $G_p$. Let us be more precise.

\begin{prop} \label{holomres} There exists a meromorphic function $\phi_p:\Si\to\End V_{\tau_p}$,
holomorphic in the region $\C_+$ if $p\not =\frac{dn}{2}$ and  in the region
$\C_+\smallsetminus\{\sqrt{\al_p}\}$ if $p=\frac{dn}{2}$, such that the resolvent $R_p(\hat s)$ is
given by the operator
\begin{eqnarray*} L^2(\Hyp^n_\K)&\longrightarrow& L^2(\Hyp^n_\K)\\ u&\longmapsto&
F_p(s,\cdot)*\phi_p(s)u
\end{eqnarray*} in the indicated regions.
\end{prop}

\begin{proof} Since the expression \eqref{eqBochner} is asymptotic to the Euclidean  one for small
$t$, we know that a radial solution of the equation
$$(\Delta_p-\alpha_p+s^2)v=0$$ must behave as $[\vol(\bS^{dn-1})t^{dn-2}]^{-1}$ as $t\to 0$ (if 
$dn>2$; the argument is similar in the other case). Thus there exists a meromorphic function
$$\psi_p\,:\,\Si\to\End V_{\tau_p},$$ holomorphic on
$\Si\priv\cN$ and such that
\begin{equation}
\label{equpper0p} F_p(\hat s,a_t)\underset{t\to 0}{\simeq}\frac{ \psi_p(\hat
s)}{\vol(\bS^{dn-1})t^{dn-2}}.
\end{equation} Consequently, for any $u\in C_0^\infty(G,\tau_p)$ we have
\begin{equation}\label{resolv} (\Delta_p-\alpha_p+s^2)(F_p(\hat s,\cdot)*u)=\psi_p(\hat s)u.
\end{equation} Moreover, for $s\in \C_+$, our previous estimates \eqref{upperp},
\eqref{equpper0p} and
\eqref{eqlowerbound0}, \eqref{eestimgreen} imply the following one:  there exists a positive
constant $C(s)$ such that
\begin{equation}\label{majoration}
\forall t>0,\ \norm{F_p(s,a_t)}_{\End V_{\tau_p}}\le C(s) G_0(\re s,a_t).
\end{equation} Hence the operator
\begin{equation}\label{borne} u\mapsto F_p(s,\cdot)*u \quad\textrm{is bounded from $L^2$ to $L^2$}.
\end{equation}

Now we study the invertibility of our function $\psi_p$.

\begin{lem}\label{petitlem}
\begin{enumerate}
\item If $p\not=\frac{dn}{2}$, the function $\psi_p$ is invertible  (with holomorphic inverse) in
the set
$\C_+$.
\item If $p=\frac{dn}{2}$, the function $\psi_p$ is invertible (with  holomorphic inverse) in the
set
$\C_+\smallsetminus\{\sqrt{\al_p}\}$.
\end{enumerate} Moreover, in both cases, $\psi_p^{-1}$ extends meromorphically to $\Si$.
\end{lem}

\begin{proof} Assume first $p\not=nd/2$. For $s\in \C_+$, let 
$\xi\in\ker \psi_p(s)$. Then
$v(a_t)=F_p(s,a_t)\xi$ provides a solution of the equation
\begin{equation}\label{equares} (\Delta_p-\alpha_p+s^2)v=0,
\end{equation} and by \eqref{estinterm} this solution satisfies
\begin{equation}\label{eqasymp} v(a_t)\underset{t\to+\infty}{=}e^{-t (\sqrt{E+s^2}+\rho 
)}\xi+\text{o}\left(e^{-t (\sqrt{E+s^2}+\rho)}\xi\right).
\end{equation} Hence $v$ is $L^2$, but we know that \eqref{equares} has no nontrivial $L^2$
solutions since
$\spec(\De_p)=[\al_p,+\infty)$ is purely continuous by \tref{specGK}. Thus
$v=0$, and therefore $\xi=0$ by \eqref{eqasymp}. It follows that
$\psi_p$ is invertible in the half-plane $\C_+$, with  holomorphic inverse in this region, and
that it has a meromorphic extension to
$\Sigma$.

Suppose now $p=\frac{dn}{2}$. Then we know that the discrete spectrum  of $\Delta_p$ reduces to
$\{0 \}$, with infinite multiplicity. Proceeding as above, we get the  second part of our lemma.
\end{proof}

According to the lemma and to \eqref{resolv}, \eqref{borne}, for 
$s\in \C_+$ (and with the additional condition $s\not=\sqrt{\al_p}$ if 
$p=\frac{dn}{2}$), the operator
$$u\in L^2\mapsto F_p(s,\cdot)*\psi_p(s)^{-1}u$$ must be the resolvent $R_p(s)$ of the operator
$\Delta_p-\alpha_p+s^2$, and this proves our proposition.
\end{proof}

We can now sum up our discussion by stating:

\begin{thm} \label{estimate} The Schwartz kernel $G_p(s,\cdot)$ of the resolvent
$$R_p(s)=(\De_p-\al_p+s^2)^{-1}$$ 
has a meromorphic extension to
$\Sigma$ and, outside the discrete subset of poles which lie inside $\Si\smallsetminus \C_+$
(except if $p=\frac{dn}{2}$, in which case $\sqrt{\al_p}\in\C_+$ is also a pole), it satisfies the
estimate
$$\forall t>1,\ \norm{G_p(\hat s,a_t)}_{\End V_{\tau_p}}\le  C(\hat s)e^{-\left(\rho +h(\hat s)
\right)t}$$ 
for some constant $C(\hat s)>0$. Moreover, we have $h(\hat s)=\re s$ on $\C_+$.
\end{thm}

\begin{rem}
1) The meromorphic extension of the resolvent $R_p$ to $\Si$ was known to 
several authors. Namely, U.~Bunke and M.~Olbrich (\cite{BO2}, Lemma~6.2) proved the
result for all hyperbolic spaces and their convex
cocompact quotients (except in the exceptional case $\K=\O$), a fact which was already observed 
by R.~Mazzeo and R.~Melrose (\cite{MM}) in the real case and by C.~Epstein, G.~Mendoza and
R.~Melrose (\cite{EMM}) in the complex case.

2) The estimate in \tref{estimate} was announced by N.~Lohou\'e in
\cite{Lohouecras}, but only in the region $\C_+$.

3) The analysis we have carried out in this section for the resolvent $R_p(s)$ is similar to the
one presented in \cite{faraut} and \cite{Miawa} for the function case ($p=0$).
\end{rem}

In order to prepare some results of next section, we discuss now the possible location of the 
poles of the resolvent $R_p$ on $\C$. We first look at the imaginary axis. 

\begin{prop} \label{poleiR}
The resolvent $R_p$ has no pole inside the set
$$i\R\smallsetminus\{\pm i\sqrt{c(\sigma_{\max})-c(\sigma)},\ \sigma\in \widehat{M}(\tau_p)\}.$$
\end{prop}

\begin{proof} Assume that $s=i\lambda$ is a purely imaginary pole of
$R_p$.
As in the proof of \lref{petitlem}, we see that there exists $\xi\in V_{\tau_p}$ such that the
function defined by $v(a_t)=F_p(i\la,a_t)\xi$ is a solution of the equation
\begin{equation}\label{eqv}
(\Delta_p-\alpha_p-\lambda^2)v=0
\end{equation}
on $\Hyp^n_\K$. 
Moreover, when $t\to+\infty$, this
solution satisfies the estimate (see \eqref{estinterm})
\begin{equation} v(a_t)\underset{t\to+\infty}{=}e^{-t (\sqrt{E-\lambda^2}+\rho
)}\xi+\text{O}(e^{-t})\xi.
\end{equation}

Let $\xi=\sum_{\sigma\in \widehat{M}(\tau_p)}  \xi_\sigma$ be the decomposition of
$\xi$ with respect to the orthogonal splitting
$$V_{\tau_p}=\bigoplus_{\sigma\in \widehat{M}(\tau_p)} V'_\sigma,
\quad\text{where}\quad V'_\si= V_\sigma\otimes\C^{m(\si,\tau_p)}.$$
Reminding \eqref{defE}, we see that there exists a certain $\eps>0$ such that the following
asymptotics holds:
\begin{equation}\label{asymptv}
v(a_t)\underset{t\to+\infty}{=}\sum_{\substack{\sigma\in \widehat{M}(\tau_p)\\ |\lambda|\ge
\sqrt{c(\sigma_{\max})-c(\sigma)}}} e^{-t (i\sqrt{\lambda^2-c(\sigma_{\max})+c(\sigma)}+\rho
)}\xi_\sigma+\text{O}(e^{-(\rho+\eps)t})\xi.
\end{equation}

Now, let $B_R$ be a geodesic ball of radius $R$ in $\Hyp^n_\K$. With the Green formula and
\eqref{eqv} we get:

\begin{eqnarray*}
0&=&\langle (\Delta_p-\alpha_p-\lambda^2)v,v\rangle_{L^2(B_R)}-
\langle v,(\Delta_p-\alpha_p-\lambda^2)v\rangle_{L^2(B_R)}\\
&=&2i\im \langle
v',v\rangle_{L^2(\partial B_R)}.
\end{eqnarray*}
But 
$$\langle v',v\rangle_{L^2(\partial
B_R)}=\vol (\partial B_R)\langle v'(a_R),v(a_R)\rangle_{V_{\tau_p}}$$ 
with 
$$\vol (\partial B_R)\simeq \vol(\bS^{dn-1}) e^{2\rho R},$$ 
so that \eqref{asymptv} implies:
$$\sum_{\substack{\sigma\in \widehat{M}(\tau_p)\\ |\lambda|\ge
\sqrt{c(\sigma_{\max})-c(\sigma)}}}
\sqrt{\lambda^2-c(\sigma_{\max})+c(\sigma)}\,\,|\xi_\sigma|^2=0.$$

Therefore, if $\lambda\not\in \{\pm \sqrt{c(\sigma_{\max})-c(\sigma)},\, \sigma\in
\widehat{M}(\tau_p)\}$, then $v \in L^2$ by \eqref{asymptv}, but we know from \tref{specGK} that
$\Delta_p$ has no $L^2$ eigenvalue inside $[\alpha_p,+\infty)$. Hence $v=0$ and this proves our
proposition.
\end{proof}

Reminding the definition \eqref{Sigma} of $\Sigma$, let us observe that the function $\hat
s\mapsto s$ is a local coordinate in a neighbourhood of
$$(0,\sqrt{c(\si_{\max})-c(\si_1)},\ldots,\sqrt{c(\si_{\max})-c(\si_r)})\in \Sigma.$$
This fact justifies the abuse of notation in the following statement.

\begin{prop} \label{zeropole} If $\al_p=0$ (i.e. if $\K=\R$ and $p=\frac{n\pm1}{2}$, see
\tref{specGK}), then the map
$s\mapsto s\,R_p(s)$ is holomorphic inside an open neighbourhood of $s=0$ in $\Si$.
\end{prop}

\begin{proof} By the spectral theorem, we know that the strong limit
$$\lim_{s\to 0^+} s^2(\Delta_p+s^2)^{-1}$$ 
is the orthogonal projector onto the $L^2$ kernel of
$\Delta_p$. But this kernel is trivial hence the limit above is zero. As we know that $s\mapsto
R_p(s)$ is meromorphic inside an open neighbourhood of $s=0$ in $\Si$, we get the result.
\end{proof}

\section{The spectrum of the differential form Laplacian on hyperbolic manifolds}

We have now all ingredients to prove the key result of our article, namely the Theorem~B
stated in the introduction, from which we shall derive various corollaries, and especially
vanishing results for the cohomology.

\subsection{Spectral results}

For convenience, let us recall here the statement of our Theorem~B. We remind
that $\de(\Ga)\le 2\rho$.

\begin{thm} \label{principal}
\begin{enumerate}
\item Assume that $p\not=\frac{dn}{2}$.
\begin{enumerate}
\item If $\delta(\Gamma)\le \rho$, then
$\lambda_0^p(\Gamma \backslash\Hyp_\K^n)\ge \alpha_p$.
\item If $\rho\le\delta(\Gamma)\le \rho+\sqrt{\al_p}$, then
$\lambda_0^p(\Gamma \backslash\Hyp_\K^n)\ge
\alpha_p-\left(\delta(\Gamma)-\rho\right)^2$.
\end{enumerate}
\item Assume that $p=\frac{dn}{2}$.
\begin{enumerate}
\item If $\delta(\Gamma)\le \rho$, then either $\lambda_0^p(\Gamma \backslash\Hyp_\K^n)=0$ or
$\lambda_0^p(\Gamma \backslash\Hyp_\K^n)\ge \alpha_p$.
\item If $\rho\le\delta(\Gamma)\le \rho+\sqrt{\al_p}$, then 
either $\lambda_0^p(\Gamma \backslash\Hyp_\K^n)=0$ or
$\lambda_0^p(\Gamma \backslash\Hyp_\K^n)\ge
\alpha_p-\left(\delta(\Gamma)-\rho\right)^2$.
\end{enumerate}
Moreover, if $\delta(\Gamma)<\rho+\sqrt{\al_p}$ the possible
eigenvalue $0$ is discrete and spectrally isolated.
\end{enumerate}
\end{thm}

When $\delta(\Gamma)>\rho+\sqrt{\al_p}$, assertions (b) are still valid, but yield a triviality
since we know that the spectrum must be non negative.

\begin{proof} Suppose first that $p\not=\frac{dn}{2}$, and let $s>0$. By our estimate
\eqref{majoration}, we have
\begin{equation}\label{aqw}
\|g_p(s,d(x,y))\|\leq C(s)\, g_0(s,d(x,y))
\end{equation}
for all $x\not=y$ in $\Hyp^n_\K$. Thus, if $s+\rho>\delta(\Gamma)$, from our \pref{prop1} we see
that, for $x\not=y$, the sum
\begin{equation}\label{resquotient}
\sum_{\gamma\in
\Gamma}\gamma_y^*\, g_p(s,d(x, y))
\end{equation}
is finite and defines therefore the
Schwartz kernel of an operator 
$$T_p(s)\,:\, C^\infty_0
\left(\wedge^pT^*( \Gamma
\backslash\Hyp_\K^n)\right)\longrightarrow
L^2_{\text{loc}}\left(\wedge^pT^*( \Gamma \backslash\Hyp_\K^n)\right).$$ 
Moreover, for any $L^2$ $p$-form $\alpha$ on $\Gamma
\backslash\Hyp_\K^n$, we have by \eqref{aqw}
$$\|T_p(s)\alpha\|_{L^2}\le C(s)\,
\|(\Delta_0-\rho^2+s^2)^{-1} u\|_{L^2}, $$
where $u=|\alpha|$. Since the operator
$\left(\Delta_0-\rho^2+s^2\right)^{-1}$ is bounded on
$L^2$, our operator $T_p(s)$ is also bounded on $L^2$. Moreover it is easy to
check that $T_p(s)$ provides a right inverse (and thus also a left inverse, by
self-adjointness) for the operator $\Delta_p-\alpha_p+s^2$ on $\Gamma
\backslash\Hyp_\K^n$. In other words, $T_p(s)$ is the resolvent of $\De_p-\alpha_p+s^2$ on 
$\Gamma \backslash\Hyp_\K^n$.

>From this discussion we see that $\la_0^p(\Gamma \backslash\Hyp_\K^n)\ge\al_p-s^2$ for any
$s>0$ such that $T_p(s)$ exists, i.e. such that $s+\rho>\de(\Ga)$. This proves assertion (1).

Suppose now that $p=\frac{dn}{2}$. The proof above still works, except when $s=\sqrt{\al_p}$, in
which case the pole of $R_p(s)$ may yield also a pole for the resolvent on the quotient. Reminding
\tref{specGK}, we get the last part of assertion (2).
\end{proof}

Before proving that a part of our estimates are optimal for a wide class of hyperbolic manifolds,
let us give some information about the nature of the spectrum of $\De_p$. Recall that
$\overline{\C_+}$ denotes the closure of $\C_+$ in $\Sigma$ (see \eqref{Sigma} and \eqref{Cplus}).

\begin{prop} \label{natspec} Assume that $\delta(\Gamma)<\rho$. The resolvent $$s\mapsto
(\Delta_p-\alpha_p+s^2)^{-1}$$ on
$\Gamma\backslash\Hyp_\K^n$, initially defined on $\C_+$, has a holomorphic
extension to an open neighbourhood of 
$$\overline{\C_+}
\priv \{\hat s\in\Sigma\,:\,s=\pm i \sqrt{c(\sigma_{\max})-c(\sigma)},\ \sigma\in
\widehat{M}(\tau_p)\}.$$  
(When $p=\frac{dn}{2}$, the value $s=\sqrt{\alpha_p}$ must be excluded also.) In particular, the
differential form spectrum of
$\Gamma\backslash\Hyp_\K^n$ is absolutely continuous on
$$[\alpha_p,+\infty) \priv \{\alpha_p+c(\sigma_{\max})-c(\sigma),\ 
\sigma\in \widehat{M}(\tau_p)\}.$$
\end{prop}

\begin{proof} When $\delta(\Gamma)<\rho$, the proof of \tref{principal} shows that the sum
\eqref{resquotient} converges for $h(\hat s)>\de(\Ga)-\rho$ as soon as the Green kernel 
$G_p(\hat s,\cdot)$ of the resolvent on $\Hyp_\K^n$ is holomorphic in the considered region.
Reminding \tref{estimate} and \pref{poleiR}, and observing that the equality \eqref{h}
extends to an open neighbourhood of 
$$\overline{\C_+}
\priv \{\hat s\in\Sigma\,:\,s=\pm i \sqrt{c(\sigma_{\max})-c(\sigma)},\ \sigma\in
\widehat{M}(\tau_p)\},$$  
we get the first assertion. The second one is 
obtained as in \cref{remcdv}.
\end{proof}

Next, let us observe that if the limit set $\Lambda(\Gamma)$ of $\Ga$ is not the whole sphere
at infinity $S^\infty=\partial \Hyp_\K^n=\bS^{dn-1}$, then
the injectivity radius of  $\Gamma\backslash\Hyp_\K^n$ is not bounded. Indeed,
if $x\in \bS^{dn-1}\priv \Lambda(\Gamma)$, then $x$ has a neighbourhood in
$\Hyp_\K^n\cup  \bS^{dn-1}$ which is 
isometrically diffeomorphic to an open subset in $\Gamma\backslash\Hyp_\K^n$ via the covering
map.

Let us remind that the condition $\Lambda(\Ga)\not=\bS^{dn-1}$ is automatically realized in
the setting of convex cocompact or geometrically finite with infinite volume quotients of
$\Hyp_\K^n$.

These remarks served us as a motivation for the two following results.

\begin{prop} If the injectivity radius of $\Gamma\backslash\Hyp_\K^n$ is not bounded (for
instance if $\Lambda(\Gamma)\not= \bS^{dn-1}$), then
$$[\alpha_p,+\infty)\subset \spec \left(\Delta_p,\Gamma\backslash\Hyp_\K^n\right)
\text{ when } p\not=\tfrac{dn}{2},$$ and
$$\{0\}\cup[\alpha_p,+\infty)\subset \spec \left(\Delta_p,\Gamma\backslash\Hyp_\K^n\right)
\text{ when } p=\tfrac{dn}{2}.$$
\end{prop}

\begin{proof}
If the injectivity radius of the Riemannian manifold $\Gamma\backslash\Hyp_\K^n$ is not
bounded, then
$\Gamma\backslash\Hyp_\K^n$ contains arbitrary large balls isometric to geodesic balls in
$\Hyp_\K^n$. But an argument due to H.~Donnelly and Ch.~Fefferman (see the proof of
Theorem~5.1.(iii) in \cite{DF}) implies then
that the essential spectrum of $\Delta_p$ on $\Gamma\backslash\Hyp_\K^n$ contains the 
essential spectrum of $\Delta_p$ on $\Hyp_\K^n$. 
\end{proof}
 
Together with \tref{principal}, this result yields immediately a generalization of a result of
R.~Mazzeo and R.~Phillips (Theorem~1.11 in \cite{MP}) when $\de(\Ga)\le\rho$.

\begin{cor} If $\delta(\Gamma)\le \rho$ and if the injectivity radius of
$\Gamma\backslash\Hyp_\K^n$ is not bounded, then 
$$\spec \left(\Delta_p,\Gamma\backslash\Hyp_\K^n\right)
=\spec \left(\Delta_p,\Hyp_\K^n\right)
=\begin{cases}
[\alpha_p,+\infty)&\text{if }p\not=\frac{dn}{2},\\
\{0\}\cup[\alpha_p,+\infty)&\text{if }p=\frac{dn}{2}.
\end{cases}$$
\end{cor}

\begin{rem} When $\Gamma$ is convex cocompact, our corollary is a particular case of a
result due to U.~Bunke and M.~Olbrich (see \S 11 in \cite{BOJFA}, and also Theorem~9.1 in
\cite{OlbrichHab}). Actually, these authors give a much more precise information: the full
spectral resolution for all vector bundles over convex cocompact quotients of $\Hyp^n_\K$ (except
for $\K=\O$ and $\de(\Gamma)\ge\rho$). 
\end{rem}

\subsection{Comparison with the Bochner-Weitzenb\"ock method}
\label{compBW}

We think it is worthwhile to compare our estimates for the bottom of the spectrum
with the ones we can get with a less elaborated method, based on the
Bochner-Weitzenb\"ock formula. Hopefully, it will turn out that the estimates in \tref{principal}
are strictly better than the latter.

Let $X=(X^m,g)$ be a complete Riemannian manifold of dimension $m$. The
\textit{Bochner-Weitzenb\"ock formula} is the identity
\begin{equation}\label{BochnerW}
\Delta_p\alpha=\nabla^*\nabla \alpha+\cour^p\alpha,\quad\forall \alpha\in
C^\infty_0(\wedge^pT^*X),
\end{equation}
where $\nabla$ is the connection on $\wedge^p T^*X$ induced by the
Levi-Civita connection and $\cour^p$ is a field of symmetric
endomorphisms of $\wedge^pT^*X$ built from the curvature tensor
(see e.g. \cite{GM}). For instance, when $X$ has
constant curvature $-1$ (typically $X=\Ga\backslash\Hyp^m_\R$), the curvature term $\cour^p$ is
quite simple:
$$\cour^p=-p(m-p){\rm Id},$$
but in general it can be hardly calculated.
Let us define thus $\cour^p_{\min}$ to be the infimum over $x \in X$ of the
lowest eigenvalues of the symmetric tensors $\cour^p(x)\,:\,
\wedge^pT_x^*X\longrightarrow \wedge^pT_x^*X$. With the Kato
inequality, we see that
$$\int_X(\Delta_p\alpha,\alpha)\ge\int_X |\nabla
\alpha|^2+\cour^p_{\min}\int_X | \alpha|^2\ge
(\lambda_0+\cour^p_{\min})\int_X | \alpha|^2$$ 
for any $\alpha\in C^\infty_0(\wedge^pT^*X)$.
In other words:

\begin{prop}\label{vanishbochner}
We have $\la_0^p(X)\ge\lambda_0^0(X)+\cour^p_{\min}$.
\end{prop}

Let us then compare the lower bounds for $\la_0^p(X)$ given by \tref{principal} and
\pref{vanishbochner}. For simplicity, we shall look only to the real and complex cases.

\subsubsection{The real hyperbolic case} 

We take $X=\Ga\backslash\Hyp^n_\R$. Recall that $\cour^p=-p(n-p)\Id$ and let us restrict to the
case $p<n/2$, thanks to Hodge duality.  By \tref{principal} and \tref{valalphap}, we have
\begin{equation*}
\la_0^p(\Ga\backslash \Hyp^n_\R)\ge
\begin{cases}
(\frac{n-1}{2}-p)^2&\text{if }\de(\Ga)\le\frac{n-1}{2},\\
(n-1-p-\de(\Ga))(\de(\Ga)-p)&\text{if }\frac{n-1}{2}\le\delta(\Gamma)\le n-1-p.
\end{cases}
\end{equation*}
Using \pref{vanishbochner} instead, we get that
\begin{equation*}
\la_0^p(\Ga\backslash \Hyp^n_\R)\ge
\begin{cases}
(\frac{n-1}{2}-p)^2-p&\text{if }\de(\Ga)\le\frac{n-1}{2},\\
\de(\Ga)(n-1-\de(\Ga))-p(n-p)&\text{if }\frac{n-1}{2}\le\delta(\Gamma)\le n-1-p.
\end{cases}
\end{equation*}
We thus see that, in all cases,
$$\text{(estimate from \tref{principal})}=\text{(estimate from \pref{vanishbochner})}+p.$$

\subsubsection{The complex hyperbolic case} 

We take now $X=\Ga\backslash \Hyp^n_\C$. In that situation, our first task is to compute the
value of $\cour^p_{\min}$:

\begin{prop} For the manifold $X=\Ga\backslash \Hyp^n_\C$, we have
$$\cour^p_{\min}=
\begin{cases}
-2p(n+1)&\text{if $p\le n$},\\-2(2n-p)(n+1)&\text{if $p\ge n$}.
\end{cases}$$
\end{prop}

\begin{proof}
In order to prove this result, we collect first some information from
the article \cite{Pedon99}.

For $0\le p\le 2n$, the representation $\tau_p$ of $K$ splits up into the direct sum 
\begin{equation*}
\tau_p=\bigoplus_{r+s=p}\tau_{r,s}
\end{equation*}
corresponding to the decomposition into differential forms of type $(r,s)$. 
Besides, each  
$\tau_{r,s}$ can be decomposed in its turn into irreducible subrepresentations: 
\begin{equation*}
\tau_{r,s}=\bigoplus_{k=0}^{\min(r,s)}\tau'_{r-k,s-k},
\end{equation*}
a fact which actually reflects the \textit{Lefschetz decomposition} into primitive forms. To sum
up, we have
\begin{equation}\label{dectaup2}
\tau_p=\bigoplus_{r+s=p}\bigoplus_{k=0}^{\min(r,s)}\tau'_{r-k,s-k}.
\end{equation}
Let us mention three natural equivalences:
\begin{align}
\tau_{r,s}&\sim\overline{\tau_{s,r}}&&\text{(complex conjugation)},\notag\\
\tau_{r,s}&\sim\tau_{n-s,n-r}&&\text{(Hodge duality)},\label{Hodgers}\\
\tau_p&\sim\tau_{2n-p}&&\text{(idem)},\label{Hodgep}
\end{align}
whose first two hold also for the $\tau'_{r,s}$.
Denoting by $\boch=\nabla^* \nabla$ the Bochner Laplacian, we can write the
following Bochner-Weitzenb{\"o}ck formulas:
\begin{align}
\De'_{r,s}&=\boch'_{r,s}+\tau'_{r,s}(\Om_\Kg)=\boch'_{r,s}-c(\tau'_{r,s})\Id,\label{W1}\\
&\qquad\qquad\text{where}\quad
c(\tau'_{r,s})=\langle\mu_{\tau'_{r,s}},\mu_{\tau'_{r,s}}+2\de_\Kg\rangle,\\
\De_{r,s}&=\boch_{r,s}+\tau_{r,s}(\Om_\Kg),\notag\\
\De_p&=\boch_p+\tau_p(\Om_\Kg).\notag
\end{align}
Our aim is thus to calculate
$$\cour^p_{\min}=\inf\{-c(\tau'_{r-k,s-k}),\ r+s=p,\,k=0,\ldots,\min(r,s)\}.$$
>From (3.8) in \cite{Pedon99} we easily see that
\begin{equation}\label{c(tau'rs)}
c(\tau'_{r,s})=2r(n-s+1)+2s(n-r+1)=2(r+s)(n+1)-4rs
\end{equation}
when $r+s\le n$, and that $c(\tau'_{r,s})$ is given by a combination of this formula with
\eqref{Hodgers} when
$r+s\ge n$. 

>From \eqref{c(tau'rs)} we deduce that:
\begin{itemize}
\item if $(r,s)$ is fixed with $r+s\le n$, then $c(\tau'_{r,s})\ge
c(\tau'_{r-k,s-k})$ for any $k\in\N$;
\item if $r+s=p\le n$ is fixed, then $c(\tau'_{r,s})=2(r+s)(n+1)-4rs$ is maximal for
$r=0$ or $s=0$, in which cases it takes the same value $2p(n+1)$.
\end{itemize}
With \eqref{dectaup2} and \eqref{Hodgep} we finally obtain the aimed result.
\end{proof}

Now we go back to our comparisons. Assume $p<n$.
By \tref{principal} and \tref{valalphap}, we have
\begin{equation*}
\la_0^p(\Ga\backslash \Hyp^n_\C)\ge
\begin{cases}
(n-p)^2&\text{if }\de(\Ga)\le n,\\
(2n-p-\de(\Ga))(\de(\Ga)-p)&\text{if }n\le\delta(\Gamma)\le 2n-p,
\end{cases}
\end{equation*}
whereas \pref{vanishbochner} yields
\begin{equation*}
\la_0^p(\Ga\backslash \Hyp^n_\C)\ge
\begin{cases}
n^2-2p(n+1)&\text{if }\de(\Ga)\le n,\\
\de(\Ga)(2n-\de(\Ga))-2p(n+1)&\text{if }n\le\delta(\Gamma)\le 2n-p.
\end{cases}
\end{equation*}
In both cases, it turns out that
$$\text{(estimate from \tref{principal})}=\text{(estimate from \pref{vanishbochner})}+p(p+2).$$

\subsection{Applications to cohomology}

We shall use in the sequel the following notation: if $X$ is any complete manifold,
\begin{align*}
H^p(X)&=\text{$p$-th de~Rham cohomology space of $X$},\\
\cH^p(X)&=\text{Hilbert space of $L^2$ harmonic $p$-forms on $X$}.
\end{align*}

Let us remark first that our \tref{principal} extends to the case of a differential form 
Laplacian
with values in a unitary flat vector bundle or even in a Hilbertian flat vector bundle. 
Indeed, if $(\pi,H_\pi)$ is a unitary representation of
$\Ga$, the $H_\pi$-valued Hodge-de~Rham operator $\De_p^\pi$ is simply
$\De_p\otimes\Id_{H_\pi}$ when lifted to the universal cover. Thus we can use similar estimates
for the corresponding Green kernel.

Keeping this generalization in mind, we state now the following vanishing result.

\begin{thm} \label{corvanish}
Assume $p\not=\frac{dn}{2}$.
If $\Ga$ and $p$ are such that $\delta(\Gamma)< \rho+\sqrt{\alpha_p}$, then 
$$\cH^p(\Gamma\backslash\Hyp_\K^n;H_\pi)=\{0\}.$$
\end{thm}

\begin{proof} When $\al_p>0$, the corollary follows immediately from \tref{principal}. When
$\al_p=0$, we know from \pref{zeropole} that  
$s\mapsto sR_p(s)=s(\De_p+s^2)^{-1}$ extends holomorphically to a neighbourhood of
$s=0$ in $\Si$. Let us denote by $\De_p^\Gamma$ the Laplacian acting on the quotient
$\Gamma\backslash\Hyp^n_\K$. Then the proof of \tref{principal} shows
that the map $s\mapsto s(\De_p^\Gamma +s^2)^{-1}$ also extends holomorphically to a neighbourhood
of $s=0$ in $\Si$ as soon as the assumption $\delta(\Gamma)<
\rho$ is fulfilled. 
As in the proof of \pref{zeropole} the spectral theorem implies that the
$L^2$ kernel of $\De_p^\Gamma$ must be trivial.
\end{proof}

\begin{rem}
1) In the convex cocompact case, this vanishing result has been also proved by M.~Olbrich
(see Corollary~9.9 in \cite{OlbrichHab}).

2) In \cite{kato}, D.~Calderbank, P.~Gauduchon and M.~Herzlich have proved refined
Kato inequalities for special classes of sections of vector bundles $E$ over a Riemannian (or
spin) manifold $X=(X^m,g)$. Namely, they consider bundles over $X$ attached to an irreducible
representation of the holonomy group $SO(m)$ and sections which lie in the kernel of a natural
injectively elliptic first-order differential operator. Their approach is based on the
representation theory of
$SO(m)$. In our situation $X=\Gamma\backslash\Hyp_\K^n$, an
application of their results (Theorem~3.1.ii and Theorem~6.3.ii) gives the following statement: 
if $\alpha$ is a harmonic $p$-form on $\Gamma\backslash\Hyp_\K^n$, then we have the
refined Kato inequality
$$|\nabla \alpha|^2\ge \frac{dn-p+1}{dn-p}\left|d| \alpha|\right|^2.$$
As a consequence, with notation of \sref{compBW}, if
$\frac{dn-p+1}{dn-p}\lambda_0^0(\Gamma\backslash\Hyp_\K^n)+\cour^p_{\min}>0$, then
$\cH^p(\Gamma\backslash\Hyp_\K^n)=\{0\}$.

When $\K=\R$ or $\C$, an easy calculation shows that this vanishing result is strictly weaker than
our \tref{corvanish}. When $\K\not=\R$, an obvious explanation is that one expects another
refined Kato inequality based on the representation theory of $K$ instead of the one of $SO(dn)$.
On the other hand, as shown in \cite{kato}, in
order to obtain an optimal result with this technique one has to consider $\Delta|\alpha|^\theta$, 
where $\alpha$ is a $L^2$ harmonic $p$-form, and
$\theta=(dn-p-1)/(dn-p)$. An easy computation shows that
$$\Delta|\alpha|^\theta\le \theta\left(-\cour^p_{\min}\right) |\alpha|^\theta.$$ 
Hence, if $|\alpha|^\theta$ is non zero and $L^2$, we get
$\lambda_0^0(\Ga\backslash\Hyp^n_\K)\le \theta\left(-\cour^p_{\min}\right)$. In that case, the
vanishing result we obtain recovers our \tref{corvanish} in the real case, and is still weaker in
the complex case. Moreover it is in general very difficult to check if $|\alpha|^\theta\in L^2$.
\end{rem}

In some cases, the Hilbert spaces $\cH^p(X)$ have a topological interpretation in terms of
cohomology groups, in the spirit of the Hodge Theorem (see for instance, \cite{Mazzeo} and
\cite{nader} for convex cocompact real hyperbolic manifolds, and \cite{MP} for geometrically finite
real hyperbolic manifolds). In this
direction, our vanishing result \tref{corvanish} also provides vanishing results for
certain cohomology groups, with a dependance on the critical exponent. 

In fact, for convex cocompact real hyperbolic manifolds, the vanishing results we can
derive from \tref{corvanish} are well known (see
\cite{izeki}, \cite{IzekiN}, \cite{Nayatani}, as well as \cite{Wang} for another approach
based on the Bochner technique). 

We therefore prefer to focus on two topological applications of our
\tref{corvanish} which seem completely new. The first one enables us to investigate the number
of ends of certain classes of hyperbolic manifolds, but we postpone its statement until next
section (see \tref{thvanish}), which will be particularly devoted to that question in a more
general setting.  The second one is specific to the complex hyperbolic case:

\begin{prop}
Suppose that $\Gamma\backslash\Hyp_\C^n$ is convex cocompact. If $p$ is such that
$p>n$ and $p>\delta(\Gamma)$, then $\cH^p(\Gamma\backslash\Hyp_\C^n)=\{0\}$
and $H^p(\Gamma\backslash\Hyp_\C^n)=\{0\}$.
\end{prop}

\begin{proof} According to T.~Ohsawa and K.~Takegoshi (Corollary~4.2  in \cite{OT}),
if $(M,h)$ is a complete Hermitian manifold of complex dimension $n$ which is K\"ahlerian
outside some compact subset $A$, and such that the K\"ahler form can be written as
$\omega=i\partial \bar\partial s$ with $s\in C^\infty(M\priv A,\R)$, $\lim_{m\to\infty}
s(m)=+\infty$ and $\partial s$  bounded, then there is an isomorphism $\cH^p(M)\simeq H^p(M)$
for any $p>n$. These assumptions may not be satisfied by $M=\Gamma\backslash\Hyp_\C^n$ with
$\Gamma$ convex cocompact. However, we claim that the complex hyperbolic metric of such
a manifold $M$ is quasi-isometric to a Hermitian metric $h$ which fulfils the above conditions;
and this is obviously enough to apply the result of T.~Ohsawa and K.~Takegoshi.

Let us elaborate. For convenience, we shall view the complex
hyperbolic space $\Hyp_\C^n$ as the open unit ball $\mathbb{B}^n_\C$ of $\C^n$. With our choice
of normalization of the Riemannian metric \eqref{bilin}, this manifold is equipped with a
K\"ahler metric of constant holomorphic sectional curvature equal to
$-4$, and the corresponding K\"ahler form is given by the formula
$$\tilde\omega=- i\partial\bar \partial \log (1-|z|^2)
=i \frac{\sum_i dz_i\wedge d\bar z_i}{1-|z|^2}
+i\frac{(\sum_i \bar z_i d z_i )\wedge (\sum_i  z_i d\bar z_i )}{(1-|z|^2)^2}.$$  
Letting $\tilde s(z)=-\log (1-|z|^2)$, we thus have
$\tilde\omega=i\partial \bar\partial \tilde s$ with $\lim_{m\to\infty} \tilde s(m)=+\infty$ and
$\partial \tilde s$  bounded. Moreover, on $\mathbb{B}^n_\C$, the (Riemannian) hyperbolic
metric
$g_{\text{hyp}}$ and the Euclidean one $g_{\text{eucl}}$ are easily compared:
\begin{equation} g_{\text{hyp}}\ge e^{\tilde s} g_{\text{eucl}}.
\label{compare}
\end{equation} 

Next, we observe that our $M=\Gamma\backslash\Hyp_\C^n$ (with $\Gamma$ convex cocompact) is
diffeomorphic to the interior of a compact manifold $\overline{M}$ with boundary $\partial
\overline{M}$, and each point $p\in
\partial \overline{M}$ has a neighbourhood $V_p$ in $\overline{M}$ which is isometric to a
neighbourhood of
$(1,0,\ldots,0)$ in $\overline{\Hyp_\C^n}=\overline{\mathbb{B}^n_\C}$. Thus, by the preceding
observation, there exists a function $s_p$ on $V_p$ such that 
$s_p^{-1}(\infty)=\partial\overline{M}\cap V_p$, $\omega=i\partial \bar\partial s_p$ and 
$\partial s_p$ is bounded (here, $\omega$ denotes the K\"ahler metric on $M$). By compactness,
we can exhibit a finite subset
$\{p_1,\ldots,p_l\}\subset\partial\overline{M}$ such that $\partial \overline{M} \subset
\bigcup_{i}V_{p_i}$. Let
$\{\varphi_i\}$ be a partition of unity associated with the covering  $\bigcup_i V_{p_i}$ and
let $s=\sum_i\varphi_i s_{p_i}$. It is clear that $\lim_{m\to \partial \overline{M}}
s(m)=+\infty$. On the other hand, each function
$\varphi_i$ is smooth on $\overline{V_{p_i}}$ and \eqref{compare} implies the estimates
$|d\varphi_i|=\mathrm{O}(e^{-s_{p_i}/2})$  and
$|\partial\bar\partial\varphi_i|=\mathrm{O}(e^{-s_{p_i}})$ on $\overline{V_{p_i}}$. Hence we
have 
$$|i\partial\bar\partial s-\omega|\le C \sum_i e^{-s_{p_i}/2} \chi_{V_{p_i}},$$ 
where $\chi_{V_{p_i}}$ denotes the
characteristic function of $V_{p_i}$. Since 
$$\lim_{m\to \partial \overline{M}} \sum_i e^{-s_{p_i}(m)/2}
\chi_{V_{p_i}}(m)=0$$ 
we find that the K\"ahler metric $\omega$ on $M$ is, near the boundary
$\partial \overline{M}$, quasi-isometric to the K\"ahler metric
$i\partial\bar\partial s$. A similar argument shows also that $\partial s$ is bounded. Thus, if
$h$ denotes a Hermitian metric on $M$ which coincides with the Hermitian metric
associated with $i\partial\bar\partial s$ near the boundary
$\partial  \overline{M}$, then $h$ is quasi-isometric to the Hermitian metric
associated with $\omega$, everywhere on $M$ (since any two Hermitian metrics are quasi-isometric
on a compact set). This discussion proves our claim.

Now, recall from \tref{valalphap} that $\alpha_p=(n-p)^2$. Since $p>n$ we have
$\delta(\Gamma)<p=n+\sqrt{\alpha_p}$, and we can apply \tref{corvanish} to obtain the
vanishing result.
\end{proof}

As a consequence, we partially recover a result of G.~Besson, G.~Courtois and S.~Gallot
(\cite{BCGacta}):
   
\begin{cor} Assume that $\Gamma\backslash\Hyp_\C^n$ is a compact complex hyperbolic manifold.
Let $\pi\,:\, \Gamma\to SU(m,1)$ be a convex cocompact representation of $\Gamma$, where $m<2n$.
Then $\delta\left(\pi\left(\Gamma\right)\right)\ge 2n=\delta\left(\Gamma\right)$.
\end{cor}

\begin{proof} From our last proposition, if
$\delta\left(\pi\left(\Gamma\right)\right)< 2n$, then
$$H^{2n}(\pi\left(\Gamma\right)\backslash\Hyp_\C^m)=\{0\}.$$ 
But, by definition of $\pi$ we have
$H^{2n}(\pi\left(\Gamma\right)\backslash\Hyp_\C^m)= H^{2n}(\Gamma\backslash\Hyp_\C^n)$, and the
latter cohomology group is obviously non trivial since
$\Gamma\backslash\Hyp_\C^n$ is a compact oriented manifold.
This discussion forces $\delta\left(\pi\left(\Gamma\right)\right)\ge 2n$. (Note that
$\de(\Ga)=2\rho=2n$ because
$\Ga$ is cocompact.)
\end{proof} 

\begin{rem}
The result of G.~Besson, G.~Courtois and S.~Gallot is in fact much better than ours: it holds
without any assumptions on $n$ and
$m$ and it says also that there is a constant $C(n,m)$ such that if
$\delta\left(\pi\left(\Gamma\right)\right)\le 2n+C(n,m)$ then
$\pi\left(\Gamma\right)$ is a totally geodesic representation. 
Note that the analogue of this phenomenon in the real hyperbolic case is also known 
(see \cite{bowen}, \cite{izeki}, \cite{IzekiN}, \cite{Nayatani}, \cite{wang}, \cite{Wang}, as well
as \cite{BCGacta} for a different proof and a more general result).
\end{rem}

\section{On the number of ends of certain noncompact locally symmetric spaces}

Let $X$ be an open manifold of dimension $m$. In what follows, we shall use the classical
notations:
\begin{align*}
H^p_0(X)&=\text{$p$-th compactly supported de~Rham cohomology space of $X$},\\
H_p(X)&=\text{$p$-th homology space of $X$}.
\end{align*}
We shall also consider the analogues of theses (co)homology spaces with coefficients in the
constant presheaf $\Z$, denoted by
$H^p_0(X,\Z)$ and $H_p(X,\Z)$, respectively. When $X$ is orientable, the Poincar\'e duality
asserts that
$$H^p_0(X)\simeq [H^{m-p}(X)]^* \simeq H_{m-p}(X),\quad H^p_0(X,\Z)\simeq
[H^{m-p}(X,\Z)]^*\simeq H_{m-p}(X,\Z),$$
as soon as these spaces are finite dimensional.

Next, recall that the \textit{number of ends} of $X$ is the
supremum over all compact subsets $A\subset X$ of the number of unbounded connected components of
$X\priv A$. 

In this section, we shall give sufficient conditions for a noncompact locally symmetric space $X$
(not necessarily of rank one) to have only one end, by showing in fact a stronger result (as is
well-known), namely that $H^1_0(X)=\{0\}$.  Our motivation was at the beginning to look
at the complex hyperbolic case, after E.~Ghys posed the problem to the first author. It
turns out that we were actually able to consider more general situations.

Before describing our results, we need some topological tools.

\subsection{Topological preliminaries}

Let us begin with the following result (see \cite{carronduke},
Theorem~3.3, for a related observation, and compare with \cite{LiWang} as well).

\begin{prop} \label{inject}
If $X=(X^m,g)$ is a complete Riemannian manifold such that every unbounded connected
component of the complement of any compact subset of $X$ has infinite volume  (for instance if
the injectivity radius is positive) and such that $\lambda_0^0(X)>0$, then the natural map
$$H^1_0(X)\longrightarrow \cH^1(X)$$ is injective. 
In particular, if furthermore $\lambda_0^1(X)>0$ then $X$ has only one end (and also 
$H_{m-1}(X)=\{0\}$ if $X$ is orientable).
\end{prop}

\begin{proof} Recall first that the spaces of $L^2$ harmonic forms admit a reduced $L^2$
cohomology interpretation:
$$\cH^p(X)\simeq \{\alpha\in L^2(\wedge^pT^*X),
d\alpha=0\}/\overline{dC^\infty_0(\wedge^{p-1}T^*X)},$$ where closure is taken with respect
to the $L^2$ topology. Hence, if $[\alpha]\in H^1_0(X)$ is mapped to zero in $\cH^1(X)$,
there is a sequence $(f_k)$ of smooth functions with compact support on $X$
such that $\alpha=\lim_{L^2} df_k$. Since we have the inequality
$$\|df_k-df_l\|_{L^2}^2\ge \lambda_0^0(X) \|f_k-f_l\|_{L^2}^2,$$ 
and since $\la_0^0(X)>0$, we conclude that this sequence $(f_k)$ 
converges to some $f\in L^2$, so that
$\alpha=df$. 
But $\alpha$ has compact support, hence $f$ is locally constant outside the compact
set $\supp(\al)$. Since all unbounded connected components of $X\priv
\supp(\al)$ have infinite volume and since $f\in L^2$, we see that $f$ has
compact support, hence $[\alpha]=[df]=0$.
\end{proof}

In the proof of the last proposition, we have used the fact that $X$ has only one end as soon as
$H^1_0(X)=\{0\}$. Next result gives a sort of converse.

\begin{prop} \label{cohomZ}
If $X=X^m$ is an open manifold having one end, and if every twofold normal covering of
$X$ has also one end, then $$H^1_0(X,\Z)=\{0\}.$$ 
In particular, $H^1_0(X)=\{0\}$ and if furthermore
$X$ is orientable, then
$$H_{m-1}(X,\Z)=\{0\}.$$
\end{prop}

\begin{proof} 
Since $X$ has only one end, we have an exact sequence
$$\{0\}\rightarrow H^1_0(X,\Z)\rightarrow H^1(X,\Z).$$ 
Pick an element in $H^1_0(X,\Z)$, and consider its image $\sigma$  in
$H^1(X,\Z)$. With $\sigma$ is associated  a continuous map $f\,: \, X\rightarrow \bS^1$, and an
induced homomorphism $f_*\,:\, \pi_1(X)\rightarrow \Z$. Because $\sigma$ has a representative with
compact support, $f$ is constant outside a compact set $C$; this constant is normalized
to be $1$.

Assume that $\sigma$ is not zero, then $f_*$ is not zero either, and has image $n\Z$, with
$n\not=0$. Then
$\Gamma=\ker\{f_* \text{ mod } 2n\Z\}$ is a normal subgroup of index
$2$ in $\pi_1(X)$. Let $\widehat X$ be the corresponding twofold normal covering of $X$, and 
let $\pi\,:\,\widehat{ X}\rightarrow X$ be the covering map. Putting $s(z)=z^2$, we have a
commutative diagram:
$$\begin{CD}
\wh{ X}@>\wh{f}>> \bS^1\\
@V{\pi}VV @VV{s}V\\
X @>{f}>> \bS^1
\end{CD}$$
But now $\widehat{X}\priv \pi^{-1}(C)$ has at least two unbounded connected
components. Indeed, on the open set $\widehat{ X}\priv \pi^{-1}(C)$,
$\widehat{f}$ is locally constant, taking both values $1$ and $-1$. Hence a contradiction,
so $\sigma$ must be trivial in $H^1(X,\Z)$, and eventually
$H^1_0(X,\Z)=\{0\}$.
\end{proof}

Although we shall not need it in the sequel, let us mention that we obtain a new proof of a
result due to Z.~Shen and C.~Sormani (\cite{SS}) as a corollary of \pref{cohomZ}.

\begin{prop}
If $X=(X^m,g)$ is a complete oriented Riemannian manifold with non negative Ricci curvature, then
either:
\begin{enumerate}
\item $H_{m-1}(X,\Z)=\{0\}$; 
\item or $X$ is the determinant line bundle of a non orientable compact manifold with
non negative Ricci curvature, and in that case $H_{m-1}(X,\Z)\simeq \Z$;
\item or $X$ is isometric to $\Sigma\times \R$ with
$\Sigma$ an oriented compact Riemannian manifold with non negative Ricci curvature, and
in that case $H_{m-1}(X,\Z)\simeq \Z$.
\end{enumerate}
\end{prop}

\begin{proof} 
According to a famous result of J.~Cheeger and  D.~Gromoll (\cite{CG}), either $X$ has one end or
$X$ is isometric to $\Sigma\times\R$, with $\Si$ as in the statement (3). 
Assume that the first possibility holds. Then we have the same alternative for any twofold
normal covering $\wh{X}$ of $X$. If $\wh X$ has only one end, we can apply \pref{cohomZ} and
obtain (1).

Thus, let us assume instead that $\wh X$ is isometric to
$\widehat{\Sigma}\times \R$, with $\wh\Si$ as before. This means
that $X=(\widehat{\Sigma}\times\R )/\{\Id,\gamma\}$ for some isometry $\ga$ of
$\widehat{\Sigma}\times \R$. By the Cheeger-Gromoll result, a line in $\widehat{\Sigma}\times
\R$ is of the form
$\{\theta\}\times \R$, where $\theta\in \widehat{\Sigma}$. Since $\gamma$ must preserve
the set of lines in $\widehat{\Sigma}\times \R$, we see that there exist $a\in\R$ and
an isometry $f$ of $\widehat{\Sigma}$ such that $\gamma(\theta,t)=(f(\theta), \pm t+a)$. 
Since also $\gamma\circ \gamma =\Id$, we must
have $\gamma(\theta,t)=(f(\theta), -t+a)$. And as $X$ is oriented, we see that $f$ has to reverse
orientation on $\widehat{\Sigma}$.
\end{proof}

\subsection{The case of general hyperbolic manifolds}

We are now able to give the second topological application of \tref{corvanish}.

\begin{thm} \label{thvanish}
Assume that $\Hyp_\K^n\not=\Hyp_\R^2$. If $\delta(\Gamma)<\rho+\sqrt{\alpha_1}$ 
and if all unbounded connected components of the complement of any compact subset of 
$\Gamma\backslash\Hyp_\K^n$ have infinite volume, then
$\Gamma\backslash\Hyp_\K^n$ has only one end, and
$$H_{dn-1}(\Gamma\backslash\Hyp_\K^n,\Z)=\{0\}.$$
\end{thm}

\begin{rem} \label{remvanish}
Except maybe for $\K=\O$, we know from \tref{valalphap} that $\al_1=(\rho-1)^2$,
hence the assumption
$\delta(\Gamma)<\rho+\sqrt{\alpha_1}$ in this statement is equivalent to 
$\delta(\Gamma)<2\rho-1$. 
Since in any case $\de(\Ga)\le 2\rho$, we see that our assumption is
not too restrictive.
\end{rem}

\begin{proof} With the hypotheses of the theorem, we know from
\tref{corvanish} that $\la_0^1(\Gamma\backslash\Hyp_\K^n)>0$. By \tref{SC}, we also
have $\la_0^0(\Gamma\backslash\Hyp_\K^n)>0$, except if $\de(\Ga)=2\rho$. Since we have assumed
$\delta(\Gamma)<\rho+\sqrt{\alpha_1}$, this cannot occur, as shown by \eqref{alphap}.

So, the result follows from \pref{inject} and \pref{cohomZ}.
\end{proof}

Actually, the assumption on $\de(\Ga)$ in the previous result is useless in the
quaternionic and octiononic cases:

\begin{cor} \label{ends}
Let $\K=\H$ or $\O$. If all unbounded connected components of the complement of any compact
subset of $\Gamma\backslash\Hyp_\K^n$ have infinite volume,
then $\Gamma\backslash\Hyp_\K^n$ has only one end, and 
$$H_{dn-1}(\Gamma\backslash\Hyp_\K^n,\Z)=\{0\}.$$
\end{cor}

\begin{proof} The hypothesis implies that
$\Gamma\backslash\Hyp_\K^n$ itself has infinite volume, so that we can use a rigidity result due to
K.~Corlette (\cite{corlette}, Theorem~4.4; see also \cite{OlbrichHab}, Corollary~4.22 for a
slight refinement).

Suppose first that $\K=\H$. 
According to Corlette's result, we have
$\delta(\Gamma)\le 4n$. In particular we always have $\delta(\Gamma)\le 4n<4n+1=2\rho-1$,
so that \tref{thvanish} applies.

Suppose now that $\K=\O$. In that case Corlette's result says that
$\delta(\Gamma)\le 16$. M.~Olbrich kindly communicated to us that he was able to calculate
the value of $\al_1$, namely he found $\alpha_1=97$, so that $\rho+\sqrt{\al_1}>20>\de(\Gamma)$ and
we can use again \tref{thvanish}. Another possible argument is the following: \tref{SC} implies
that $\lambda_0^0(\Gamma\backslash\Hyp_\O^2)\ge 96=6\times 16$, and since
$\Gamma\backslash\Hyp_\O^2$ is an Einstein manifold with Ricci curvature equal to $-36$,
the Bochner formula \eqref{vanishbochner} yields $\la_0^1(\Gamma\backslash\Hyp_\O^2)\ge 60>0$,
so that $\cH^1(\Gamma\backslash\Hyp_\O^2)=\{0\}$. Thus we can use \pref{inject} and
\pref{cohomZ}.
\end{proof}

\begin{rem}
Our \cref{ends} extends a result of K.~Corlette about convex cocompact quotients of quaternionic
and octonionic hyperbolic spaces (see \cite{corlette}, Theorem~7.1).
\end{rem}

As another consequence of \tref{thvanish}, we give a simple proof of a result due to
Y.~Shalom (\cite{shalom}, Theorem~1.6), which we shall actually improve
a bit later on in the $SU(n,1)$ case (see \cref{sha2}).

\begin{cor} \label{sha1}
Assume that $\Ga=A*_C B$ is a cocompact subgroup of $SO_e(n,1)$ (with $n\ge 3$) or $SU(n,1)$
(with $n\ge 2$) which is a free product of subgroups $A$ and $B$ over the amalgamated subgroup $C$.
Then $\delta(C)\ge 2\rho-1$.
\end{cor}

\begin{proof} Let $\K$ be either $\R$ or $\C$.
By a recent result of G.~Besson, G.~Courtois and S.~Gallot (\cite{BCGpre}) we have
$$H_{dn-1}(C\backslash\Hyp_\K^n)\not=\{0\}.$$ 
But $C\backslash\Hyp_\K^n$ is a Riemannian
covering of a compact hyperbolic manifold, so its injectivity radius has a uniform positive
lower bound and the unbounded connected components of the complement of any
compact subset of
$C\backslash\Hyp_\K^n$ must have infinite volume. To avoid contradiction with  
\tref{thvanish} (and Remark~\ref{remvanish}), we must have $\delta(C)\ge 2\rho-1$.
\end{proof}

\begin{rem} In his paper, Y.~Shalom proves actually a better
result in the complex case, namely, that the inequality is
strict. Besides, \cite{BCGpre} gives a substantial generalization of Shalom's result: if
$A*_C B$ is the fundamental group of a compact Riemannian manifold $(X^m,g)$ with sectional
curvature less than $-1$, then $\delta(C)\ge m-2$. Also, the equality case is characterized
when $X$ is real hyperbolic and $m\ge 4$.

Note that \cref{sha1} is meaningless in the quaternionic or octonionic case.
Indeed, since $Sp(n,1)$ and $F_{4(-20)}$ satisfy the property (T) of Kazhdan, it is
well known that none of their cocompact subgroups can be an almagamated product (see \S 6.a
in \cite{T}).
\end{rem}

\subsection{The case of locally symmetric spaces which have the Kazhdan property}

Let us give now an analogue of \cref{ends} in the case of more general noncompact locally
symmetric spaces whose isometry group satisfies Kazhdan's property $(T)$.

\begin{thm} \label{endT}
Let $G/K$ be a symmetric space without any compact factor and without any factor
isometric to a real or complex hyperbolic space. Assume that $\Gamma\subset G$ is a torsion-free,
discrete subgroup of $G$ such that $\Gamma\backslash G/K$ is non compact and that all unbounded
connected components of the complement of any compact subset of
$\Gamma\backslash G/K$ have infinite volume. Then
$\Gamma\backslash G/K$ has only one end, and $$H_{m-1}(\Gamma\backslash G/K,\Z)=\{0\},$$
where $m=\dim(G/K)$.
\end{thm} 

\begin{proof} Under our assumptions $G$ satisfies property $(T)$, and the
quotients $\Gamma\backslash G/K$ and
$\Gamma\backslash G$ have infinite volume. Thus the right regular representation of $G$ on
$L^2(\Gamma\backslash G)$ has no nontrivial almost invariant vector, and this implies that
$\lambda_0^0(\Gamma\backslash G/K)>0$: if instead we had
$\lambda_0^0(\Gamma\backslash G/K)=0$, we could construct a sequence $(f_l)$ of smooth functions
with compact support on
$\Gamma\backslash G/K$ such that $\|df_l\|_{L^2}\le\|f_l\|_{L^2} /l$. By pulling back this
sequence to $\Gamma\backslash G$, we would obtain a sequence of nontrivial
almost invariant vectors in $L^2(\Gamma\backslash G)$, which is absurd.

Next, the fact that $\cH^1(\Gamma\backslash G /K)=\{0\}$, and thus that $\Gamma\backslash G
/K$ has only one end by \pref{inject}, is also a heritage of the property $(T)$. 
Let us elaborate. 

According to N.~Mok (\cite{Mok}) and P.~Pansu (\cite{Pansu}), the property $(T)$ for the
group $G$ can be shown with a Bochner type formula which is in fact a special case of a
refinement of the Matsushima formula obtained by N.~Mok, Y.~Siu and S.~Yeung 
(\cite{mys}). In particular there exists on
$G/K$ (and on $\Gamma\backslash G/K$) a parallel curvature tensor $B$ which is positive definite
on symmetric $2$-tensors having vanishing trace, and such that for any 
$L^2$ harmonic $1$-form $\al$ on $\Gamma\backslash G/K$ we have:
\begin{equation}\label{B}
\int_{\Gamma\backslash G/K} B(\nabla \alpha,\nabla \alpha)d\vol =0.
\end{equation}
Since $\alpha$ is closed and coclosed, $\nabla \alpha$ is symmetric and has vanishing trace,
thus formula \eqref{B} implies that $\alpha=0$; hence $\cH^1(\Gamma\backslash G /K)=\{0\}$. 
Note that \eqref{B} is usually stated in the finite volume setting.
But the extension to noncompact $\Gamma\backslash G/K$ presents no difficulties: if $\alpha$ is
a $L^2$ harmonic $1$-form on $\Gamma\backslash G/K$, it is easy to check that $\nabla\alpha$
is also $L^2$; thus, the
integration by part procedure required to derive \eqref{B} can be justified by standard cut-off
arguments.

Since our discussion clearly applies to any finite covering of $\Gamma\backslash G/K$, we finish
the proof by employing \pref{cohomZ}.
\end{proof}

\begin{rem} 
Assume instead that $\Gamma\backslash G/K$ has finite volume. If we have also $\rank_{\mathbb{Q}}
\Gamma\ge 2$, the Borel-Serre compactification of $\Gamma\backslash G/K$ implies that
$\Gamma\backslash G/K$ has only one end.
\end{rem}

\subsection{The specific case of complex hyperbolic manifolds}
Besides the result of \tref{ends}, we have for complex hyperbolic manifolds the following
statement.

\begin{thm}\label{endC}
Let $\Gamma$ be a discrete and torsion-free subgroup of $SU(n,1)$, with $n\ge 2$. Assume that the
limit set
$\Lambda(\Gamma)$ is not the whole sphere at infinity
$\bS^{2n-1}$, that $\delta(\Gamma)<2n$, and that the injectivity radius of
$\Gamma\backslash\Hyp_\C^n$ has a positive lower bound. Then
$\Gamma\backslash\Hyp_\C^n$ has only one end, and
$$H_{2n-1}(\Gamma\backslash\Hyp_\C^n,\Z)=\{0\}.$$
\end{thm}

Note that the hypotheses in this theorem are always satisfied in the
convex cocompact setting.

\begin{proof} By \tref{SC}, the hypothesis $\delta (\Gamma)<2n$ implies that 
$c=\la_0^0(\Gamma\backslash\Hyp_\C^n)>0$. Thus the following Poincar\'e inequality holds:
\begin{equation}\label{poin}
\forall f\in C^\infty_0(\Gamma\backslash\Hyp_\C^n),\quad c\|f\|^2_{L^2}\le \|df\|^2_{L^2}.
\end{equation}
On the other hand, our assumption on the injectivity radius 
implies that the volume of geodesic balls of radius $1$ is uniformly bounded from below.
Since the Ricci curvature of $\Gamma\backslash\Hyp_\C^n$ is constant, a result by N.~Varopoulos 
(see \cite{varo}, or Theorem~3.14 in \cite{Hebey}) asserts that, for some other constant $c'>0$,
we have the Sobolev inequality:
\begin{equation}\label{sob}
\forall f\in C^\infty_0(\Gamma\backslash\Hyp_\C^n),\quad c'\|f\|^2_{L^{n/(n-1)}}\le
\|df\|^2_{L^2}+\|f\|^2_{L^2}.
\end{equation} 
Gathering inequalities \eqref{poin} and \eqref{sob}, we obtain the following Euclidian type
Sobolev inequality: for some constant $c''>0$,
\begin{equation}\label{sobeu}
\forall f\in C^\infty_0(\Gamma\backslash\Hyp_\C^n),
\quad c''\|f\|^2_{L^{n/(n-1)}}\le \|df\|^2_{L^2}.
\end{equation}

Next, suppose that there exists a compact set $C\subset \Gamma\backslash\Hyp_\C^n$ such that
$(\Gamma\backslash\Hyp_\C^n)\priv C$ has at least two unbounded connected components, and let
$\Omega$ be one of them. 
According to Theorem~2 in \cite{CSZ}, thanks to
\eqref{sobeu} we can find a harmonic function $u$ on
$\Gamma\backslash\Hyp_\C^n$, which is valued in $[0,1]$ and satisfies
$$\int_{\Gamma\backslash\Hyp_\C^n} |du|^2d\vol<+\infty,$$ 
as well as 
\begin{equation}\label{limu}
\lim_{\substack{m\to\infty\\ m\in \Omega}}u(m)=0
\ \text{ and } \lim_{\substack{m\to\infty\\ m\not\in \Omega}} u(m)=1.
\end{equation} 
By Lemma~3.1 in \cite{Li}, $u$ must be pluriharmonic. In particular, $u$ is
harmonic on any complex submanifold of $\Gamma\backslash\Hyp_\C^n$.

Now, let $p\in\bS^{2n-1}\priv\Lambda(\Gamma)$. Then there exists a neighbourhood $U$ of $p$ in
$\Hyp_\C^n\cup\bS^{2n-1}$, such that $U$ is mapped isometrically in
$\Gamma\backslash\Hyp_\C^n$ by the covering map 
$\pi\,:\,\Hyp_\C^n\rightarrow \Gamma\backslash\Hyp_\C^n$. 
But we can find a holomorphic map $F\,:\,{\mathbb
D}\rightarrow U$ such that $F(\partial{\mathbb D})=F({\mathbb D})\cap
\bS^{2n-1}$. For instance, if $p=(1,0,...,0)$ then for some $\varepsilon >0$ small enough,
$$z\mapsto F(z)=(\sqrt{1-\varepsilon^2},\varepsilon z,0,...,0)$$ 
is such a map. So $u\circ\pi\circ F$ is
a bounded  harmonic function on ${\mathbb D}$, and takes a constant value on $\partial{\mathbb D}$
($0$ or $1$). Hence $u$ is constant on 
$\pi\circ F({\mathbb D})$ and, by the Maximum Modulus Theorem, $u$ must be constant everywhere.
This contradicts \eqref{limu}, so that $\Gamma\backslash\Hyp_\C^n$ must have only one end.

The vanishing result follows again from \pref{cohomZ}.
\end{proof}

\begin{rem} 
Actually the proof of \tref{endC} extends to the case of any complete K\"ahler manifold containing
a proper holomorphic disc and verifying the Sobolev estimate \eqref{sobeu}.
We recover thus a result of J.~Kohn and H.~Rossi (\cite{kohn}) which
asserts that a K\"ahler manifold which is pseudo-convex at infinity has only one end. 
There is a lot of literature which deals with the number of ends of complete K\"ahler
manifolds, see for instance the references \cite{LiR} and \cite{NR}.
\end{rem}

As an immediate consequence of our last theorem, we can complement the result of Y.~Shalom that we
recovered in \cref{sha1}:

\begin{cor} \label{sha2}
Assume that $\Gamma=A*_C B$ is a cocompact subgroup of $SU(n,1)$ (with $n\ge 2$) which is a free
product of subgroups $A$ and $B$ over an amalgamated subgroup $C$. Then either
$2n-1\le\delta(C)<2n$ and $\Lambda(C)=\bS^{2n-1}$, or $\delta(C)=2n$.
\end{cor}


\end{document}